\documentclass[twoside,leqno]{article}

\usepackage[letterpaper]{geometry}

\usepackage{siamproceedings}

\usepackage[T1]{fontenc}
\usepackage{amsfonts}
\usepackage{times} 
\usepackage{graphicx}
\usepackage{epstopdf}
\usepackage{enumitem}
\usepackage{algorithmic}
\ifpdf
  \DeclareGraphicsExtensions{.eps,.pdf,.png,.jpg}
\else
  \DeclareGraphicsExtensions{.eps}
\fi

\newsiamremark{remark}{Remark}
\newsiamremark{hypothesis}{Hypothesis}
\crefname{hypothesis}{Hypothesis}{Hypotheses}
\newsiamthm{claim}{Claim}

\usepackage{amsopn}
\DeclareMathOperator{\diag}{diag}

\usepackage{amssymb}
\usepackage{mathrsfs}
\usepackage{xcolor}

\global\long\def\bx{\mathbf{x}}%
\global\long\def\bxs{\mathbf{x}_{\star}}%
\global\long\def\bF{\mathbf{F}}%
\global\long\def\by{\mathbf{y}}%
\global\long\def\bv{\mathbf{v}}%
\global\long\def\bu{\mathbf{u}}%
\global\long\def\gradsp{\nabla_{\rm sp}}%
\global\long\def\Hesssp{\nabla^{2}_{\rm sp}}%
\global\long\def\qs{q_{\star}}%
\global\long\def\Es{E_{\star}}%
\global\long\def\BN{\mathbb{B}^{N}}%
\global\long\def\SN{\mathbb{S}^{N-1}}%
\global\long\def\E{\mathbb{E}}%
\global\long\def\P{\mathbb{P}}%
\global\long\def\R{\mathbb{R}}%
\global\long\def\indic{\mathbf{1}}%
\global\long\def\TAP{{\rm TAP}}%
\global\long\def\crt{{\rm Crt}}%
\global\long\def\B{{\rm Band}}%

\begin{document}

\newcommand\relatedversion{}
\renewcommand\relatedversion{\thanks{The full version of the paper can be accessed at \protect\url{https://arxiv.org/abs/0000.00000}}} 

\title{\Large Geometry of spherical spin glasses\relatedversion}
    \author{Eliran Subag\thanks{Department of Mathematics, Weizmann Institute of Science (\email{eliran.subag@weizmann.ac.il}).\\
    	\hspace*{.52cm}The author was partially supported by ERC grant PolySpin 101165541 and ISF grant 2055/21.}
    }

\date{}

\maketitle


\fancyfoot[R]{\scriptsize{Copyright \textcopyright\ 20XX by SIAM\\
Unauthorized reproduction of this article is prohibited}}





\begin{abstract} 
	Spherical spin glasses are canonical models for smooth random functions in high dimensions. In this review, we survey several interrelated lines of research on their geometric structure. We begin with results concerning critical points and their relationship to the Gibbs measure. For the pure models, the measure concentrates on spherical bands around critical points that approximately maximize the energy at a particular radius. 
	Next, we present another approach in which a similar picture is derived for general mixed models. At the core of this approach is a free energy functional computed over bands using multiple orthogonal replicas, satisfying a strong concentration of measure. We discuss several implications of this method for a generalized Thouless-Anderson-Plamer (TAP) approach.
	Finally, we explain how these geometric insights inform optimization algorithms, and briefly relate them to Smale’s 17th problem over the real numbers.
\end{abstract}

\section{Introduction.}

Magnetic materials are classified into several types based on their response to external magnetic fields.
The most well-known are ferromagnetic materials, such as iron, in which spins of neighboring atoms tend to align in the same direction.  
In contrast, in antiferromagnetic materials, spins align in alternating directions and cancel each other.
Other types include ferrimagnetic, paramagnetic and diamagnetic materials, each with unique spin behavior.

In the early 1970s, experimental physicists began to observe anomalous magnetic behavior in certain disordered alloys, which did not fit any known category. These materials were  named \emph{spin glasses} by analogy with the positional disorder of atoms in structural glasses. Their properties could not be explained by existing theories, and they introduced a new paradigm in the study of disordered systems.
The discovery raised a remarkable new problem for theoretical physicists. 
To start making sense of the experiments and develop a theory, they needed new models. 

In 1975, Edwards and Anderson (EA) \cite{EA1975} proposed a modification of the well-known Ising model for ferromagnets and antiferromagnets where interaction coefficients have a fixed sign, by making them random to account for the irregular, distance-dependent interactions between atoms of the alloy. In the same year, Sherrington and Kirkpatrick (SK) \cite{SK75} introduced a mean-field model of a spin glass, where any two spins interact with each other, unlike the EA model in which interactions are only between neighboring sites on a grid. Although less realistic, the mean-field SK model was expected to exhibit the same qualitative behavior as in the EA model (though this remains debated to this day).

We will focus on the \textit{spherical mixed $p$-spin models}. To choose a specific model from this class we fix a function  $\xi(t):=\sum_{p=1}^{\infty}\gamma_{p}^{2}t^{p}$ called the \emph{mixture}, where for simplicity we assume that $\xi(t)<\infty$ for some $t>1$. The corresponding spherical mixed $p$-spin \emph{Hamiltonian}, or \emph{energy}, is the zero-mean Gaussian process on the sphere $$\SN:=\left\{\bx \in\mathbb{R}^{N}:\|\bx\|_{2}=\sqrt{N}\right\}$$ with covariance function
\begin{equation}
	\E\left\{ H_{N}(\bx)H_{N}(\by)\right\} =N\xi(R(\bx,\by)),\label{eq:cov}
\end{equation}
where $R(\bx,\by):=\bx\cdot\by/N$ is  the  \emph{overlap
	function.} More explicitly, it can be realized as
\begin{equation}
	H_{N}(\bx)=\sum_{p=1}^{\infty}\gamma_{p}H_{N,p}\left(\bx\right):=\sum_{p=1}^{\infty}\gamma_{p}N^{-\frac{p-1}{2}}\sum_{i_{1},...,i_{p}=1}^{N}J_{i_{1},...,i_{p}}^{(p)}x_{i_{1}}\cdots x_{i_{p}}\,,\label{eq:mixed}
\end{equation}
where $\bx=(x_{1},\ldots,x_{N})$ and the disorder coefficients $J_{i_{1},...,i_{p}}^{(p)}$
are i.i.d. standard normal variables.
For $\xi(t)=t^p$ the model is called the \emph{pure} $p$-spin model. More classically, mixed models were defined on the hyper-cube $\Sigma_N=\{-1,1\}^{N}$ in which case the spins are said to be Ising spins.  The SK model mentioned above is the pure $2$-spin model with Ising spins. 
	
	Their importance in physics is not the only reason to study spin glass models (although the principle that natural problems in physics often lead to elegant, interesting mathematics has definitely proved itself in their case). In this review paper we take the perspective that they are canonical models for the study of smooth random functions in high dimensions.
	To elaborate on this point, we make two observations. The first, and obvious, observation is that Hamiltonians as in \eqref{eq:mixed} are random  \emph{polynomials} if $\gamma_p=0$ for $p>p_0$. The set of all such polynomials is of course dense in the space of continuous functions on $\SN$. 
		The second remark is that \eqref{eq:cov} means that $H_N(\bx)$ is a sequence of stationary processes on the sphere $\SN$, with  covariance having the same dependence on the normalized distance at all dimensions. 
In fact, by a classical result of Schoenberg \cite{Schoenberg}, mixtures of the form above are the most general functions $\xi(t)$ such that the right-hand side of \eqref{eq:cov} is a covariance function (i.e., positive-definite for all $N$)\,---\,at least up to adding a $p=0$ term, which is immaterial to us as it amounts to shifting the whole process by the same Gaussian variable.

Consider the level sets 
\begin{equation}
	L_{N}(E):=\left\{ \bx\in\SN:\, H_{N}(\bx) = NE\right\} \label{eq:lvl}
\end{equation}
as a tool to investigate the geometry of the Hamiltonian. 
They are not empty with high probability (w.h.p.) for $E\in (-\Es,\Es)$, where $\Es$, called the \emph{ground state} energy, is the almost sure (a.s.) non-random limit 
\begin{equation}
	\Es:=\lim_{N\to\infty}\frac{1}{N}\max_{\bx\in \SN}H_{N}(\bx)\,.\label{eq:GS}
\end{equation}

As random subsets of a high-dimensional manifold with no obvious special structure, $L_N(E)$ are extremely complicated. We start with a coarse yet very challenging geometric question: what is the asymptotic behavior of their volumes $\text{Vol}(L_{N}(E))$? 
 It is immediate to see that  they decay exponentially with $N$ (unless $E=0$).\footnote{For super-level sets, exponential decay of the normalized $N-1$ dimensional volume follows by Markov's inequality and Fubini. For the level sets, exponential decay of the $N-2$ dimensional volume follows from a variant of the Kac-Rice formula.} The right scale is therefore to consider
\begin{equation}
	V_N(E):=\frac1N\log \Big(\text{Vol}(L_{N}(E))\Big)\label{eq:vol}.
\end{equation}
In fact, and this is highly non-trivial to prove, (see Footnote \ref{ft:3}) 
the almost sure limit 
$
V(E):=\lim\limits_{N\to\infty}V_N(E)
$ 
exists. 

It is sometimes useful to approximate volumes by integrals of smooth functions. For reasons that will gradually become clearer, we work with exponential functions and analogously to $V_N(E)$ and $V(E)$ define
\begin{equation}
	F_{N,\beta}:=\frac{1}{N}\log Z_{N,\beta}:=\frac{1}{N}\log\bigg(\int_{\SN}e^{\beta H_{N}(\bx)}d\bx\bigg),\label{eq:FN}
\end{equation}
where the integration is w.r.t. the uniform probability measure on the sphere, and
\begin{equation*}
	F(\beta):=\lim_{N\to\infty}F_{N,\beta}\,.
\end{equation*}
These are classical objects in statistical physics. Namely, $Z_{N,\beta}$ and $F_{N,\beta}$ are called the \emph{partition function} and \emph{free energy}, respectively, while the parameter $\beta\geq0$ represents the \emph{inverse temperature}.

By the coarea formula, neglecting the polynomial in $N$ Jacobian, we have that
\[
Z_{N,\beta} \approx \int_{-\infty}^{\infty} e^{\beta NE} \text{Vol}(L_{N}(E))dE = \int_{-\infty}^{\infty} e^{N(\beta E+V_N(E))} dE\,.
\]
Actually, one can show that
\[
F(\beta) = \max_{|E|<\Es}\Big(\beta E + V(E)\Big)\qquad\mbox{and}\qquad V(E)=\min_{\beta\geq0}\Big(F(\beta)-\beta E\Big)\,.
\]
In other words, the limiting free energies $F(\beta)$ and volumes $V(E)$ contain the same information (as functions).

Perhaps even more interesting than the volume is the structure of $L_N(E)$, which one may investigate through the uniform measure on it. 
By analogy to the above, we consider a related `smoothed' measure called  the \emph{Gibbs measure}
\begin{equation}
	G_{N,\beta}(A):=\frac{1}{Z_{N,\beta}}\int_{A}e^{\beta H_{N}(\bx)}d\bx\,,\qquad A\subset \SN\,.\label{eq:GN}
\end{equation}
If we choose $\beta$ such that $E=F'(\beta)$ then, for some  $\varepsilon_N\to0$, w.h.p. (see \cite{AuffingerChenConcentration} by Auffinger and Chen)
\[
G_{N,\beta}\Big(
\big\{\bx\in\SN: |N^{-1}H_N(\bx)-E|<\varepsilon_N\big\}
\Big)\approx1\,.
\]
Hence, if the level set does not change drastically as the energy varies slightly around $E$, then $G_{N,\beta}$ should capture the behavior of the uniform measure on the level set of $E$.

In 1979, Parisi discovered his celebrated formula \cite{ParisiFormula,Parisi1980} for the free energy of the Sherrington-Kirkpatrick model. In the following years, considerable effort was devoted to reformulating his solution\,---\,based on the non-rigorous \emph{replica trick}\,---\,in terms of the system's physical properties. This effort gave rise to a rich and intricate theory that offered a detailed description for the structure of  the Gibbs measure \cite{MPVspinglass}. Another very influential paper \cite{TAP} was written by Thouless, Anderson and Palmer in 1977. The method introduced there, known as the TAP approach, will be central to Section \ref{sec:TAP}. Over the past few decades, an active area of research in probability theory has emerged around the mathematical study of spin glasses, and many predictions from physics have been rigorously proved \cite{PanchenkoBook,TalagrandOldBook,TalagrandBookI,TalagrandBookII}. Most notably, Talagrand established a rigorous proof of Parisi’s formula \cite{Talag,Talag2}, while Panchenko resolved the ultrametricity conjecture \cite{ultramet}. A more detailed discussion of these developments is deferred to Section~\ref{sec:background}.

Since the dimension of the configuration space varies, in mean-field spin glass models a thermodynamic limit cannot be defined in the usual sense. Instead, a different concept of an asymptotic Gibbs measure has emerged. Precisely, let $\bx^1,\bx^2,\ldots$ be an infinite sequence of independent samples from $G_{N,\beta}$ and consider the infinite overlap array $Q_N=(R(\bx^i,\bx^j))_{i,j\geq 1}$. Under a technical condition \eqref{eq:generic}, the overlap array $Q_N$ weakly converges to an explicit random limit $Q_\infty$ characterized by Ruelle probability cascades \cite{Ruelle}\,---\,random hierarchical probability measures defined on a separable Hilbert space\,---\,see Section \ref{sec:background} for more details. 

The overlap array $Q_N$ determines the Gibbs measure $G_{N,\beta}$  up to orthogonal transformations (i.e., modulo rotations), see \cite[Lemma 1.7]{PanchenkoBook}. However, the asymptotic Gibbs measure\,---\,defined on an abstract Hilbert space\,---\,does not retain precise information about where in space the finite-$N$ Gibbs measure concentrates. A central focus of the results reviewed in this paper is precisely this issue: describing the regions where $G_{N,\beta}$ concentrates in terms of simple statistics of the Hamiltonian $H_N(\bx)$, defined on the finite-dimensional configuration space $\SN$ before the limit is taken.


The first such description we discuss is given in Section \ref{sec:crtpts} about critical points of $H_N(\bx)$. Define the spherical band 
\begin{equation}\label{eq:band}
	\B(m,\delta) := \big\{\bx\in\SN:\,|R(\bx,m)-R(m,m)|\leq \delta\big\}\,.
\end{equation}
For small $\delta>0$, it consists of points in $\SN$ that are approximately orthogonal to $m$.
For the pure models and mixed models sufficiently close to pure, if $\beta$ is large, the Gibbs measure concentrates on such thin spherical bands with $\delta=\delta_N\to0$,  centered around the few highest critical points of $H_N(\bx)$ on $\sqrt{q(\beta)}\SN = \big\{\bx \in\mathbb{R}^{N}:\|\bx\| =\sqrt{q(\beta)N}\big\}$. These bands are the `pure-states' of the Gibbs measure and their centers are approximately orthogonal to each other as 1-step replica symmetry breaking occurs (this terminology will be explained below).  
The proof methods of these results\,---\,which rely on the Kac-Rice formula and various estimates\,---\,fail in more general settings.

Another, more general, approach which covers all mixtures $\xi(t)$ is discussed in Section \ref{sec:TAP}. An interesting question is to characterize the 
points $m$ inside the sphere $\|m\|<\sqrt N$ for which $G_{N,\beta}(\B(m,\delta))$ does not decay exponentially. If an additional condition on how the Gibbs measure is distributed over the band is imposed, there is a surprisingly simple characterization: the centers of bands $m$ lie in $\sqrt{q}\SN$ for specific radii and approximately maximize $H_N(\bx)$ over the same sphere. These radii correspond to values $q\in[0,1]$ that exhibits a property we call \emph{multi-samplability}, roughly meaning that the probability that many independent samples from $G_{N,\beta}$ have pairwise overlap $q$ is not exponentially small. After presenting this characterization, we show how it naturally gives rise to a generalized TAP correction and representation for the free energy and discuss how these results relate to the Parisi solution.

One of the main appeals of spin glass theory is its broad influence across many other fields.\footnote{See \cite{MPVspinglass,Anderson1989} for an early perspective on this point, or the recent Parisi's 2021 Nobel Prize press release, which connects his 1980s contributions
to ``very different areas, such as mathematics, biology,
neuroscience and machine learning''.} Of special
importance is its connection to non-convex optimization, established with the discovery of simulated annealing
by Kirkpatrick, Gelatt and Vecchi \cite{Kirkpatrick1983} in 1983, and highlighted in the classical M\'{e}zard-Parisi-Virasoro book \cite{MPVspinglass} (where one
of the three parts is devoted entirely to optimization) as well as the series of short articles by Anderson \cite{Anderson1989}.
In Section \ref{sec:algorithms} we explain how results from the previous paragraph naturally lead to algorithms for optimization of spin glass Hamiltonians. Recent related results about Smale's 17th problem \cite{SmaleProblems} over the reals are also discussed briefly there.

\section{Additional Background: the Parisi solution and the Gibbs measure.\label{sec:background}}
The replica trick in physics suggests  computing the free energy using the identity\footnote{The first equality simply uses the identity $\log x=\lim\limits_{n\to0}(x^n-1)/n$ for $x>0$, the second follows from dominated convergence since $(x^n-1)/n$ increases with $n$ and the third follows from a Taylor expansion of $\log(x)$ around $1$.}
$$\E\log Z_{N,\beta} = \E \lim_{n \to 0}\frac1n(Z_{N,\beta}^n-1) = \lim_{n \to 0}  \frac1n(\E Z_{N,\beta}^n-1) =\lim_{n \to 0}\frac1n \log \E Z_{N,\beta}^n\,.$$
Loosely speaking, the basic idea is to write a formula for \[\lim\limits_{N\to\infty}\frac{1}{nN} \log \E Z_{N,\beta}^n\] for \emph{integer} $n\geq 1$, 
extend it to real $n$ in a way that allows  taking the limit as $n\to0$, and declare this limit as the free energy. To justify the crucial extension to the reals, analytic continuation is commonly invoked but not substantiated rigorously \cite{MPVspinglass}.
Using Fubini's theorem, the integer moments can be written as integrals over $n$ points in space, and their logarithmic $N\to\infty$ limit can be expressed by a maximum over 
$n\times n$ matrices corresponding to their overlaps.

An attempt to compute the free energy of the SK model was made by Sherrington and Kirkpatrick \cite{SK75} using the replica trick with a replica symmetry Ansatz, meaning that the aforementioned maximization is restricted to constant matrices.
As they pointed out themselves, their solution could not be true at low temperature (large $\beta$) as the entropy became negative. Another issue of stability was raised in \cite{Almeida1978}, and their result did not agree with computer simulations for the ground state energy \cite{KirkpatrickPhysRevB.17.4384}. Further attempts by Blandin \cite{Blandin} and Bray and Moore \cite{BrayMoorePhysRevLett.41.1068} made progress but did not completely resolve those issues. Finally, in 1979, Parisi invented his celebrated replica symmetry breaking (RSB) scheme \cite{ParisiFormula,Parisi1980}, which in particular identified the structure of matrices one should maximize over. 
His ideas were unconventional even by the standards of physicists at the time and mathematically very far from rigorous. However, no flaws were found in his solution, which showed excellent agreement with computer simulations and was quickly accepted by the physics community \cite{CharbonneauReplicaReview}.
Among his many remarkable contributions that earned him the 2021 Nobel Prize in Physics, his RSB solution, along with his other work on spin glasses, truly stand out \cite{ParisiNobelLecture}.

The \emph{Parisi formula} states that
\begin{equation}\label{eq:Parisi}
	F(\beta) = \inf_{\mu} 	\mathcal{P}_{\xi,\beta}(\mu)\,,
\end{equation}
where the infimum is over all probability measures supported on $[0,\hat q]$ for some $\hat q<1$. For the spherical models, $\mathcal{P}_{\xi,\beta}(\mu)$ admits the \emph{Crisanti-Sommers
	representation} \cite{Crisanti1992} (for the Ising case, see the works cited below)
\begin{equation}
	\mathcal{P}_{\xi,\beta}(\mu)=\frac{1}{2}\Big(\beta^{2}\int_{0}^{1}\mu([0,q])\xi'(q)dq+\int_{0}^{\hat{q}}\frac{dq}{\int_{q}^{1}\mu([0,s])ds}+\log(1-\hat{q})\Big)\label{eq:crisantisommers}\,.
\end{equation}
The unique minimizer in \eqref{eq:Parisi} is called the \emph{Parisi measure}. If it consists of a single atom, $k+1$ atoms, or has infinite support we say that the model exhibits replica symmetry (RS), $k$-step replica symmetry breaking ($k$-RSB), or infinite replica symmetry breaking ($\infty$-RSB). 

Rigorously, even the existence of the limit of the free energy $F(\beta)$ was not known for many years, until it was proved in 2002 by Guerra and Toninelli \cite{GuerraToninelli} (see also \cite{FreeEnergyConvergence}).\footnote{\label{ft:3}In addition to the existence of $F(\beta)$, its differentiability was proved in the Ising case in \cite{MR2195333,MR2399285}, while it is trivial for the spherical models using the same argument. Combining these with large deviations techniques as in \cite{AuffingerChenConcentration} and basic Gaussian processes theory, one can show that the $N\to\infty$ limits $V(E)$ and $\Es$ as in \eqref{eq:GS} exist.} In a major breakthrough in the same year, Guerra \cite{GuerraBound} invented an ingenious interpolation technique and proved that Parisi's formula is an upper bound to the free energy. Within less than a year, Talagrand proved the matching lower bound in his celebrated works \cite{Talag,Talag2}  and completed the proof of the formula for even mixtures $\xi(t)=\xi(-t)$.  
 The later proof of Parisi's formula for general mixtures 
$\xi(t)$, due to Panchenko \cite{Panchenko2014} and Chen \cite{Chen2013} (for the Ising and spherical models, respectively), built upon two major advances.
The first was the representation of the free energy established in the breakthrough work of Aizenman, Sims, and Starr \cite{AS2scheme}. The second was Panchenko’s celebrated resolution of the Parisi ultrametricity conjecture \cite{ultramet}, which was crucial for understanding the structure of the Gibbs measure.

Overlaps of `replicas'\,---\,i.i.d. samples $\bx^1,\bx^2,\ldots$ from the Gibbs measure $G_{N,\beta}$\,---\,have turned out to play a central role. 
Parisi  realized already in \cite{ParisiOrderPar}  that the annealed overlap distribution 
$	\mu(\cdot) = \P(R(\bx^1,\bx^2)\in\cdot )$
should coincide with the minimizer of \eqref{eq:Parisi}. Assuming the Ansatz from \cite{ParisiFormula,Parisi1980}, this implicitly suggested that w.h.p.
\begin{equation}\label{eq:ultR}
	R(\bx^1,\bx^2)\geq R(\bx^1,\bx^3) \wedge R(\bx^2,\bx^3) - o(1)\,,
\end{equation} a relation that came to be known as the Parisi ultrametricity conjecture mentioned above.  In the seminal works \cite{MPSTV2,MPSTV1},  M\'{e}zard, Parisi, Sourlas, Toulouse  and   Virasoro  explicitly put forward ultrametricity \eqref{eq:ultR} as a key structural feature in the organization of pure states and explored its implications. 

To simplify the discussion, assume that $\xi(t)$ is generic, namely
\begin{equation}\label{eq:generic}
	\sum_{p\text{ even}}p^{-1}\indic\{\gamma_{p}>0\}=\sum_{p\text{ odd}}p^{-1}\indic\{\gamma_{p}>0\}=\infty\,,
\end{equation}
in which case the infinite overlap array $Q_N=(R(\bx^i,\bx^j))_{i,j=1}^{\infty}$ weakly converges to a random array $Q=(Q_{ij})_{i,j=1}^\infty$, see \cite[Section 3.7]{PanchenkoBook}.\footnote{The proof of such convergence actually relies on the results we review next, but as we mentioned this assumption is made for simplification. Without it, in what follows one may discuss subsequential limits of $Q_N$.} The limiting array is clearly positive-definite and exchangeable\,---\,that is, invariant in law under permutations of the indices\,---\,and therefore admits a Dovbysh-Sudakov representation \cite{Dovbysh-Sudakov}. In the limit, the Parisi ultrametricity conjecture \eqref{eq:ultR} becomes the statement that almost surely
\begin{equation}\label{eq:ultQ}
Q_{12}\geq Q_{13}\wedge Q_{23}\,.
\end{equation}
Two important, closely related stability properties were established for the limiting infinite overlap array $Q$ in the late 1990s: the Aizenman-Contucci stochastic stability \cite{AizenmanContucci} and the Ghirlanda-Guerra identities \cite{GhirlandaGuerra}. Baffioni and Rosati \cite{BaffioniRosati} showed that if the infinite overlap array $Q$ satisfies the Ghirlanda-Guerra identities and is ultrametric \eqref{eq:ultQ}, then the annealed overlap distribution $\P (Q_{12}\in\cdot )$ uniquely determines
the law of $Q$.
Another crucial ingredient in the description of the Gibbs measure are the probability cascades introduced by Ruelle \cite{Ruelle}. Their infinite overlap array is ultrametric and satisfies the Ghirlanda-Guerra identities, and thus any array with both these properties either is equal in distribution to the array of the Ruelle probability cascade with the same annealed overlap distribution or is a limit of such (in case that the support of this distribution is an infinite set).
The last piece of the puzzle was Panchenko's celebrated proof \cite{ultramet} that the Ghirlanda-Guerra identities, in fact, imply ultrametricity. Consequently, the limiting overlap array $Q$ is also ultrametric and its law is characterized by Ruelle probability cascades. 

For a comprehensive treatment of the results discussed above\,---\,and much more\,---\,see the excellent books  by Talagrand \cite{TalagrandBookI,TalagrandBookII} and Panchenko \cite{PanchenkoBook}. For the interested reader, we also recommend  their 1998 and 2022 ICM papers on the state of the art nearly three decades ago \cite{MR1648045} and on ultrametricity in spin glasses \cite{MR4680407}.

\section{Critical points: complexity and the Gibbs measure.\label{sec:crtpts}} 
The focus of this section is critical points of the Hamiltonian $H_N(\bx)$ on $\SN$. Subsection \ref{subsec:complexity} deals with counts of critical points at a given energy, which grow exponentially with $N$. In Subsection \ref{subsec:Gibbsgeometry} we discuss structure theorems for the Gibbs measure at large $\beta$, identifying the critical points and regions around them at which $G_{N,\beta}$ concentrates.

\subsection{Annealed and quenched complexity and the extremal process.\label{subsec:complexity}} 

Let $F_{1},...,F_{N-1}$ be a piecewise smooth frame field on $\SN$ and define the spherical gradient and Hessian,
\begin{align*}
	\gradsp H_{N}(\bx) :=\left( F_{i}H_{N}(\bx)\right)_{i\leq N-1}, \qquad	\Hesssp H_{N}(\bx  :=\left( F_{i}F_{j}H_{N}(\bx)\right)_{i,j\leq N-1}\,.
\end{align*}
For $B\subset\mathbb{R}$, we denote the number of critical points with normalized energy in $B$ by
\begin{equation*}
	\crt_N(B):=\#\left\{ \bx\in\SN:\,\gradsp H_{N}(\bx)=0,\,H_{N}(\bx)\in NB
	\right\}\,. \label{eq:0204-01}
\end{equation*}

In the seminal works of Auffinger, Ben Arous and \v{C}ern\'{y} \cite{A-BA-C} about pure models and Auffinger and Ben Arous \cite{ABA2} about mixed models, they computed the expected complexity of critical points with energy in $B$ for any spherical model $\xi(t)$. For the explicit expression for the complexity function  $\Theta_{\xi}(u)$ below, see Eq. (1.22) in \cite{ABA2}.

\begin{theorem}[Annealed complexity, Auffinger-Ben Arous \cite{ABA2}]\label{thm:1stmoment}
	For any interval $B\subset\R$,
	\begin{equation}\label{eq:complexity}
	\lim_{N\to\infty}\frac1N\log \E\crt_N (B) = \sup_{u\in B}\Theta_{\xi}(u)\,.
	\end{equation}
\end{theorem}

In fact, finite index complexity functions were also calculated in \cite{ABA2,A-BA-C} (i.e., asymptotics for critical points with $\Hesssp H_{N}(\bx)$ having a specified number of negative eigenvalues). Calculations similar to the above were done by Fyodorov for spherical models in \cite{Fyodorov2013} and models in $\R^N$ in \cite{Fyodorov2004} and with Williams in \cite{MR2363390}.

The calculation of the complexity \eqref{eq:complexity}  in \cite{ABA2,A-BA-C} starts by applying the well-known Kac-Rice formula \cite{RFG,Kac1943,Rice1945} to express the mean by a conditional expectation involving $H_N(\bx)$, $\gradsp H_{N}(\bx)$ and the determinant of $\Hesssp H_{N}(\bx)$ with $\bx\in\SN$ arbitrary, due to stationarity.  As the latter variables are jointly Gaussian, a straightforward covariance calculation gives their joint law, see \cite[Lemma 1]{ABA2} and \cite[Section 4.1]{geometryMixed}. This results in transferring the problem of computing the mean in \eqref{eq:complexity} to a computation involving a modified GOE matrix (scaled and shifted by a random multiple of the identity matrix) the asymptotics of which are studied in \cite{ABA2,A-BA-C} using tools from random matrix and large deviations theory.

Define $E_0 = E_0(\xi)$ as the maximal value $E$ such that $\Theta_\xi(E)\geq0$ and note that $\Es\leq E_0$, by Markov's inequality.
However, the annealed complexity of Theorem \ref{thm:1stmoment} does not necessarily describe the typical behavior of $\crt_N (B)$.  Interestingly, for some explicit function $G(x,y)$, it was proved in \cite[Proposition 6]{ABA2} that $G(\xi'(1),\xi''(1))\geq0$ if and only if $E_0$ is equal to the ground-state energy obtained from the Parisi formula \eqref{eq:Parisi} restricted to 1-RSB measures for large $\beta$. More precisely, if one takes the minimum only over measures $\mu$ consisting of two atoms. Specifically, for the 1-RSB pure $p$-spin models $\Es=E_0$, as already proved in \cite[Theorem 2.12]{A-BA-C} using the Parisi formula.

In fact, for the pure models, concentration around the mean was proved for $\crt_N(B)$ by a different method. Define $E_\infty = E_\infty(p) := 2\sqrt{(p-1)/p} < E_0$.

\begin{theorem}[Quenched complexity \cite{2nd}]\label{thm:2ndmoment}
	For the pure $p$-spin model with any $p\geq3$ and $E\in(E_\infty,E_0)$,
	\begin{equation}\label{eq:momentmatch}
	\lim_{N\to\infty}\frac{\E \big\{ (\crt_N((E,\infty)))^2  \big\} }{\big( \E \{\crt_N((E,\infty))\}\big)^2 }= 1\,.
	\end{equation}
	Consequently, in $L^2$ and in probability,
	\[
	\lim_{N\to\infty}\frac{\crt_N((E,\infty))  }{ \E \{\crt_N((E,\infty))\} }= 1\,.
	\]
\end{theorem}
The same result was proved by Zeitouni and the author in \cite{2ndarbitraryenergy} for $p\geq32$ and any $E\in(-\infty,E_\infty)$, while Ben Arous, Zeitouni and the author proved in \cite{geometryMixed} matching of moments on exponential scale (i.e., that $\frac1N\log$ of the ratio in \eqref{eq:momentmatch} goes to $0$) for mixed models which are close to pure in a suitable sense and $E$ close to $E_0$. We note that by the Borell-TIS inequality, such exponential matching implies that $\Es=E_0$.

As for the first moment, the calculation of the second moment starts by applying the Kac-Rice formula, using that $(\crt_N((E,\infty)))^2$ is the number of pairs $(\bx,\by)\in \SN\times\SN$ of critical points. However, the resulting random matrix problem in this case is more challenging. It involves two modified GOE matrices whose correlations depend on an overlap parameter $r\in(-1,1)$. Asymptotic analysis of this problem \cite{2nd,geometryMixed} yields, for general mixtures $\xi(t)$,  that
\begin{equation}\label{eq:2nd_formula}
\lim_{N\to\infty} \frac1N \log(\E \crt_N(B)^2) = \sup_{r\in[0,1]}f_{\xi,B}(r)
\end{equation}
for an explicit function $f_{\xi,B}(r)$. Proving matching of moments on exponential scale, a for specific models $\xi(t)$, boils down to showing $f_{\xi,B}(r)$ attains its maximum at $r=0$.

Formulae similar to \eqref{thm:1stmoment} and \eqref{eq:2nd_formula} can be derived for critical points with the radial derivative restricted to a subset $D$ instead of the energy $H_N(\bx)/N$  \cite{geometryMixed}. Interestingly, for the analogue of \eqref{eq:2nd_formula} for this problem, Belius and Schmidt \cite{BeliusSchmidt} showed that maximality at $r=0$  holds for any mixture $\xi(t)$ if the minimum of $D$ is large enough, implying matching of moments on exponential scale.
Given a spherical model, the ground-state energy $\Es$ and the normalized radial derivative at the ground-state configuration, which we denote by $R_\star$, can be computed from the Parisi formula. In a beautiful paper, Haung and Sellke \cite{HuangSellke2023} proved that for 1-RSB models, the annealed complexity of critical points approximately with energy $\Es$ and radial derivative $R_\star$ is zero. By a modified second moment argument they also prove concentration and thus identify the value of the ground-state energy.
They use this with the generalized TAP representation we review in Section \ref{sec:TAP}, a decomposition of $\xi(t)$ to sub-models from \cite{FElandscape} and the calculation of the ground-state of full-RSB in \cite{GSFollowing}, to prove the lower bound of Parisi's formula.
Additional first and second moment computations were carried out using the Kac-Rice formula in several closely related settings including the spiked tensor model \cite{MR4011861,PhysRevX.9.011003}, topological triviality \cite{MR4404771,MR4354698,MR3162549,Huang_toptriv},
bipartite and multispecies spherical models \cite{MR3249903,DartoisMcKenna,MR4645713,MR4718393}, random autonomous ODEs \cite{MR4305622,MR3598251,MR3521630,Garciacomplexity,MR4824712,MR4411381}, the elastic manifold \cite{MR4552227,MR4673883}, Ising TAP complexity \cite{BeliusTAP1,MR4203332},   empirical risk minimization \cite{pmlr-v107-maillard20a,Montanari_emprisk}, random polynomial systems \cite{Azais2005,Subag2023,Wschebor2005}.

We conclude by discussing critical values close to the maximal energy, which will be particularly important in the next subsection. They can be described through the extremal point process of critical values
\[
\zeta_N := \iota_p^{-1} \sum_{\bx\in\SN:\,\gradsp H_N(\bx)=0}\delta_{H_N(\bx)-m_N}\,.
\]
Here, assuming that $\xi(t)=t^p$ is the pure $p$-spin model, $\iota_p\in\{1,2\}$ with the same parity as $p$ to account for the fact that critical values come in antipodal pairs in the even case, and $m_N=NE_0-\frac12\log N/c_p-K_p$ where $c_p=|\Theta_\xi'(E_0)|$ and $K_p$ is an explicit constant (see Eq. (2.6) of \cite{pspinext}). Zeitouni and the author proved the following weak convergence, where ${\rm PPP}(\mu)$ denotes the law of a Poisson point process of intensity $\mu$.

\begin{theorem}[Limiting extremal process, Subag-Zeitouni \cite{pspinext}]\label{thm:extremal}
	For the pure $p$-spin model with any $p\geq3$,
	\begin{equation*}
		\zeta_N\overset{N\to\infty}{\longrightarrow}\zeta_\infty\sim {\rm PPP}(e^{-c_px}dx)\,.
	\end{equation*}
\end{theorem}

%

\subsection{Decomposition of the Gibbs measure around critical points.\label{subsec:Gibbsgeometry}}
Recall the definition of $\B(m,\delta)$ in \eqref{eq:band} and let $\bxs^i$ be an enumeration of the critical points of $H_N(\bx)$ on $\SN$ such that $H_N(\bxs^i)\geq H_N(\bxs^{i+1})$. In Eq. (56) of \cite{geometryGibbs} 
a certain overlap value $q_*=q_*(p,\beta)$ was defined. Here we work with a slightly different parametrization and use $\qs  = \qs(p,\beta) :=q_*^2$. 
The following is the main result of \cite{geometryGibbs}.
\begin{theorem}
	\label{thm:Geometry} Suppose $\xi(t)=t^p$ is the pure $p$-spin model with odd $p\geq3$. For large enough
	$\beta$, we have that 
		\begin{equation}
			\lim_{k,c\to\infty}\liminf_{N\to\infty}\E\left\{ G_{N,\beta}\left(\cup_{i\leq k}\B\big(\sqrt{\qs}\bxs^i,cN^{-1/2}\big)\right)\right\} =1\label{eq:thm1_1}
		\end{equation}
		and for any $i,j\geq1$ and $c,\,\delta>0$
		\begin{equation}
			\lim_{N\to\infty}\E\left\{ G_{N,\beta}^{\otimes 2}\Big( \big|R(\bx^1,\bx^2)-\delta_{ij}\qs \big|>\delta\,\Big|\,\bx^1\in \B\big(\sqrt{\qs}\bxs^i,cN^{-1/2}\big),\bx^2\in \B\big(\sqrt{\qs}\bxs^j,cN^{-1/2}\big) \Big) \right\} =0.
			\label{eq:1402}
		\end{equation}	
\end{theorem}
Here $\delta_{ij}$ is the Kronecker delta. In fact, for $i=j$ a stronger bound was proved in \cite{geometryGibbs} with $\delta$ replaced by any sequence $\delta_N$ such that $N^{1/2}\delta_N\to\infty$. 
The theorem is stated for odd $p$ for simplicity, but it also holds (including the latter improvement) for even $p\geq 4$ with obvious modifications to account for the fact that since $H_N(\bx)=H_N(-\bx)$ is even, critical points and the bands around them come in identical pairs related by reflection. We also remark that the critical points that are relevant to the decomposition above are the ones characterized in Theorem \ref{thm:extremal} about the extremal process.

The  decomposition of Theorem \ref{thm:Geometry} gives a precise description of the pure states \cite{TalagrandPstates,JagannathApxUlt} of the Gibbs measure, which we discuss in Subsection \ref{sec:tree}. However, unlike the more abstract decomposition of \cite{TalagrandPstates,JagannathApxUlt}, here we give an explicit description of  the states in terms of the highest critical points.
Of course, $\qs\to1$ as $\beta\to \infty$ and the bands shrink to the critical points $\bx^i$, which are approximately orthogonal in $\R^N$.
A related technique was used by Ben Arous and Jagannath \cite{BenArousJagannathShattering} to prove shattering for $G_{N,\beta}$ with $\beta<\beta_c$.

Let $T(A,B)$ be the set of points $\bx\in\SN$ such that $R(\bx,\bxs)\in B$ for some critical point $\bxs\in\SN$ with  
$\frac1NH_N(\bxs)\in A$, and consider the contribution of this
subset to the partition function
\begin{equation}
	Z(A,B) = \int_{T(A,B)} e^{\beta H_N(\bx)}d\bx\,.
\end{equation} 
If for any $\delta>0$ and $A_\delta=E_0+[-\delta,\delta]$ and $B_\delta=\sqrt{\qs}+[-\delta,\delta]$ we have that $Z(A_\delta,B_\delta)$ is asymptotically much greater than the same integral over $\SN\setminus T(A_\delta,B_\delta)$, then 
\begin{equation}\label{eq:asymptoticsupport}
	\lim_{\delta\to0}\liminf_{N\to\infty}\E\left\{ G_{N,\beta}\left(\cup\,\B\big(\sqrt{\qs}\bxs,\delta \big)\right)\right\} =1 \,,
\end{equation} 
where the union is over critical points with $\frac1NH_N(\bxs)\in A_\delta$. 

We now briefly discuss some of the ideas in \cite{geometryGibbs} that went into the proof of this coarse version of \eqref{eq:thm1_1}. First, a lower bound on $Z(A_\delta,B_\delta)$ is proved.
It particularly implies that $F(\beta)> \beta \bar E$, for a suitable energy $\bar E=\bar E(\beta)$.
Let $L_N^+(E)$ and $L_N^-(E)$ be the super-level and sub-level sets of $H_N(\bx)$ on $\SN$ defined by
\[L_N^\pm(E) = \{\bx\in\SN: \pm H_N(\bx)\geq \pm NE\}\,.\]
Then, since the free energy concentrates \cite[Theorem 1.2]{PanchenkoBook} and 
\[
\int_{L_N^-(\bar E)} e^{\beta H_N(\bx)}d\bx \leq e^{N\beta \bar E} =  Z_{N,\beta} \exp\Big\{ N(\beta\bar E-F_{N,\beta})\Big\} \,,
\]
to prove \eqref{eq:asymptoticsupport} it is enough to bound the same integral over $L_N^+(\bar E)\setminus T(A_\delta,B_\delta)$.

For sufficiently large $\beta$, $\bar E$ is close to the ground-state energy $E_0=\Es$ and all critical points in $L_N^+(\bar E)$ are local maxima \cite{ABA2}. In fact, it was shown in \cite{geometryGibbs} that for an appropriate $\bar q\in(0,1)$, w.h.p. for any critical point $\bxs\in L_N^+(\bar E)$,
\[
\sup_{\bx\in \B(\bar q^{1/2}\bxs,0)}H_N(\bx)<N\bar E\,.
\]
Thus, each connected component of the super-level set $L_N^+(\bar E)$ contains exactly one critical $\bxs$ point and is contained in the set of points $\bx$ such that $R(\bx,\bxs)\geq \bar q^{1/2}$. In particular, $L_N^+(\bar E)\setminus T(A_\delta,B_\delta)$ can be covered by a union of sets of the form $T(A,B)$ for $A\subset [\bar E,E_0+\delta]$ and $B\subset [\bar q,1]$. Importantly, for $A$ and $B$ in this range, it can be shown that for large $\beta$, $Z(A,B)$ is an additive functional in $A$ and $B$ (while for general $A$, $B$ this is not true). The proof of \eqref{eq:asymptoticsupport} is completed by proving upper bounds for $Z(A,B)$ corresponding to a partition of $L_N^+(\bar E)\setminus T(A_\delta,B_\delta)$.

Crucially, since it is defined through critical points, $Z(A,B)$ could be analyzed using the Kac-Rice formula. If $A_N$ and $B_N$ shrink to $E$ and $\sqrt q$ as $N\to\infty$, $T(A_N,B_N)$ is a collection of thin bands, each can be approximated by a sphere $\B(q^{1/2}\bxs,0)$ of co-dimension $1$. Their analysis via the Kac-Rice formula requires estimating the free energy of a spherical Hamiltonian with a mixture that depends on $q$ (see \eqref{eq:xiq}).

Our discussion so far in this subsection has concerned the pure $p$-spin models. To generalize the results to the mixed case, it has turned out that one should center the bands around critical points of $H_N(\bx)$ on $\sqrt{\qs}\SN$ for a suitable $\qs\in(0,1)$. In the special case of the pure (homogeneous) models, this coincides with our discussion above, since such critical points correspond to those on $\SN$ simply by scaling. 
Ben Arous, Zeitouni and the author \cite{geometryMixed} proved the following for mixed models close to pure. Let $\qs=\qs(\beta)\in(0,1)$ be the square of the value defined in Eq. (8,6) of \cite{geometryMixed} and let $\bxs^i(q)$ be an enumeration of the critical points of $H_N(\bx)$ on $\sqrt q\SN$ such that $H_N(\bxs^i(q))\geq H_N(\bxs^{i+1}(q))$. Recall the definition of $E_0(\xi)$ as the zero of the annealed complexity and note that the restriction of $H_N(\bx)$ to $\sqrt{q}\SN$ corresponds to the mixture $\xi(qt)$.
For a mixture $\nu(t)=\sum_{p=2}^{\infty}\gamma_p^2 t^p$, define $\|\nu\|:=\sum_{p=2}^{\infty}\gamma_p^2 p^4$. 
\begin{theorem}[Ben Arous-Subag-Zeitouni \cite{geometryMixed}]\label{thm:Geometrymixed}
	For any $p\geq3$ there exists some $\epsilon>0$ such that if the mixture $\xi(t)=\sum_{k\geq2}\gamma_k^2t^k$ is not even and satisfies $\|\xi(t)-t^p\|<\epsilon$, then
	for large enough
	$\beta$ and some $\epsilon_N,\delta_N\to0$ we have
	\begin{align*}
		\lim_{N\to\infty}\E\max_{i\leq \exp(N\epsilon_N)}\Big|\frac1NH_N(\bxs^i(\qs))-E_0(\xi(\qs t))\Big| &= 0\,,\\
		\lim_{N\to\infty}\E\left\{ G_{N,\beta}\left(\cup_{i\leq \exp(N\epsilon_N)}\B\big(\bxs^i(\qs),\delta_N\big)\right)\right\} &=1\,,
	\end{align*}
and for any $\delta>0$,
\begin{equation*}
	\lim_{N\to\infty}\E\bigg\{ \sum_{i,j\leq \exp(N\epsilon_N)} G_{N,\beta}^{\otimes 2}\Big( \big|R(\bx^1,\bx^2)-\delta_{ij}\qs \big|>\delta,\,\bx^1\in \B\big(\bxs^i(\qs),\delta_N\big),\bx^2\in \B\big(\bxs^j(\qs),\delta_N\big) \Big) \bigg\} =0.
\end{equation*}

\end{theorem}

\section{Free energies on bands and a generalized TAP approach.\label{sec:TAP}} For the convenience of the reader, we recall 
\begin{equation}\label{eq:band2}
	\B(m,\delta) := \big\{\bx\in\SN:\,|R(\bx,m)-R(m,m)|\leq \delta\big\}\,.
\end{equation}
In Theorems \ref{thm:Geometry} and \ref{thm:Geometrymixed} we saw that for some models, the Gibbs measure concentrates on such bands centered at critical points of approximately the ground-state energy. Unfortunately, the methods of proof do not generalize to arbitrary mixtures. In this section we discuss a more general approach, which allows us to establish a similar decomposition (and much more) for arbitrary mixtures $\xi(t)$.

To motivate the discussion, we start with a natural problem of characterizing the points $m\in\BN:=\{\bx:\,\|\bx\|\leq\sqrt{N}\}$ whose band  has significant Gibbs mass for $\delta\ll 1$. Restricting our discussion to the exponential scale (i.e., the level of free energies), the question is: can we identify all points $m$ such that 
\begin{equation}
	\label{eq:Gband}
	F_{N,\beta}	\approx  F_{N,\beta}(m,\delta):=\frac1N\log\int_{\B(m,\delta)}e^{\beta H_N(\bx)}d\bx\,?
\end{equation}

One potential issue is that the set of such points $m$ may be too large to be described by any succinct or structurally meaningful condition.
For example, imagine a hypothetical scenario in which $G_{N,\beta}$ concentrates on a few points, and any band containing at least one of these points is heavy. Another serious issue is that even if we wish to analyze $F_{N,\beta}(m,\delta)$ only up to small errors of order $O(1)$, we basically need to control exponentially many points in $\BN$ and a union bound would not suffice (since $F_{N,\beta}(m,\delta)$ has exponential tails). 

Observe that the bands in Theorems \ref{thm:Geometry} and \ref{thm:Geometrymixed} satisfy an additional property beyond \eqref{eq:Gband}: two independent samples from the band typically have overlap $R(\bx^1,\bx^2)$ approximately equal to the self-overlap $R(m,m)$ of the center point $m$. In fact, the same holds for pure states \cite{TalagrandPstates,JagannathApxUlt} with $m$ the barycenter under $G_{N,\beta}$. Moreover, the same occurs for any pair among $k\geq2$ independent samples. On exponential scale, this property is that for small $\epsilon$ and large $k$, 
\begin{equation}
	\label{eq:orthpts}
	\frac1N\log G_{N,\beta}^{\otimes k}\Big(\forall i\neq j:\,
	\,\big|R(\bx^i,\bx^j)-R(m,m)\big|<\epsilon \,\Big|\,\forall i:\,\bx^i\in \B(m,\delta)
	\Big)\approx 0\,.
\end{equation}
Let us revise our previous question and now ask:  can we identify all points $m$ such that both \eqref{eq:Gband} and  \eqref{eq:orthpts} hold for $\delta,\epsilon\ll1$ and $k\gg 1$? The answer is yes, and we will see two related characterizations. First, these are points $m\in\BN$  such that
\begin{equation}
	\label{eq:char}
	\frac{\beta }{N}H_N(m)+F_{\TAP}(\|m\|^2/N)\approx F(\beta),
\end{equation}
for a deterministic function $F_{\TAP}=F_{\TAP,\beta}(q)$, see Subsection \ref{sec:TAPcorrection}. Second, they are points such that 
\begin{equation}
	\label{eq:char2}
	\|m\|\approx \sqrt{Nq}\mbox{\qquad and \qquad}
	\frac{1 }{N}H_N(m)\approx \frac{1 }{N}\max_{\|\bx\|=\sqrt{Nq}}H_N(\bx),
\end{equation}
for some \emph{multi-samplable} overlap $q\in[0,1]$, which we define in Subsection \ref{subsec:multisamplable}. 

\subsection{\label{sec:many_orth_reps}Using many orthogonal replicas.} 
Define the free energy
\begin{equation}
	\label{eq:TAPk}
	\begin{aligned}
		F_{N,\beta,k}(m,\delta,\epsilon)&:=\frac{1}{kN}\log\int_{\B_k(m,\delta,\epsilon)}e^{\beta\sum_{i\leq k}[H_N(\bx^i)-H_N(m)] }
		d\bx^1\cdots d\bx^k\,,\\
		\B_k(m,\delta,\epsilon)&:=\Big\{
		(\bx^1,\ldots,\bx^k)\in \B(m,\delta): \forall i\neq j,\, |R(\bx^i,\bx^j)-R(m,m)|<\epsilon
		\Big\}\,.
	\end{aligned}
\end{equation}
For small $\delta$ and $\epsilon$, note that  $(\bx^i-m)_{i\leq k}$ as above are approximately orthogonal to each other and to $m$. In this notation, 
\begin{equation}\label{eq:FbandF1}
F_{N,\beta}(m,\delta) = \frac{\beta}{N} H_N(m)+F_{N,\beta,1}(m,\delta,\epsilon)	
\end{equation}
(independently of $\epsilon$) and the left-hand side of \eqref{eq:orthpts} is equal to 
\begin{equation}\label{eq:Gk}
	k\big(F_{N,\beta,k}(m,\delta,\epsilon)-F_{N,\beta,1}(m,\delta,\epsilon)
	\big)\,.
\end{equation} 
Therefore, for any $m$,
\begin{equation}\label{eq:FkF1inequality}
	\frac{\beta }{N}H_N(m)+F_{N,\beta,k}(m,\delta,\epsilon)\leq \frac{\beta }{N}H_N(m)+F_{N,\beta,1}(m,\delta,\epsilon)\leq F_{N,\beta}
\end{equation} 
while both conditions \eqref{eq:Gband} and \eqref{eq:orthpts} hold if and only if
\begin{equation}\label{eq:FkF1apx}
	\frac{\beta }{N}H_N(m)+F_{N,\beta,k}(m,\delta,\epsilon)\approx\frac{\beta }{N}H_N(m)+F_{N,\beta,1}(m,\delta,\epsilon)\approx F_{N,\beta}\,.
\end{equation}

A crucial property of the free energy defined above is the following concentration result.
\begin{theorem}[Uniform concentration, {\cite[Proposition 1]{FElandscape}}] \label{thm:concentration}For any fixed  $\beta,c,t>0$, if $\delta,\epsilon>0$ are small enough and $k\geq1$ is large enough, then for large $N$,
	\begin{equation}\label{eq:uniformconcentration}
		\P\bigg(
		\max_{\|m\|\leq\sqrt N}\Big|\,
		F_{N,\beta,k}(m,\delta,\epsilon) - \E F_{N,\beta,k}(m,\delta,\epsilon)
		\,\Big|<t
		\bigg)>1-e^{-cN}.
	\end{equation}
\end{theorem}

Denote by $D(m)$ the deviation as in \eqref{eq:uniformconcentration} without the maximum. By standard concentration results \cite[Theorem 1.2]{PanchenkoBook}, $\P(D(m)\geq t)\leq 2\exp(-t^2k^2N/4a)$ where $Na$ is the maximal variance of $\beta\sum_{i\leq k}(H_N(\bx^i)-H_N(m))$ over $\B_k(m,\delta,\epsilon)$. As can be easily verified, $a/k^2$ is as small as we wish, provided that $\epsilon,\delta$ are small enough and  $k$ is large enough (uniformly in $m$). That is, for a single point we can get exponential concentration with a constant as large as we wish. Of course, this is sufficient for controlling deviations over an exponentially large net. Uniform concentration as above, is proved by showing in addition a Lipschitz property for $D(m)$, with high probability.

\subsection{Multi-samplable overlaps.\label{subsec:multisamplable}} We are interested in identifying points $m$ that satisfy \eqref{eq:FkF1apx} or, equivalently, both \eqref{eq:Gband} and \eqref{eq:orthpts}.  Thanks to the inequalities \eqref{eq:FkF1inequality}, in the limit of small $\epsilon,\delta$, large $k$ and large $N$, such points are found at radii $\sqrt{Nq}$ for
$q\in[0,1)$ such that
\begin{equation}\label{eq:eta}
\inf_{\epsilon,\delta,\eta,k}\lim_{N\to\infty}\P\Big(
\exists m\in \sqrt{q}\SN:\,\Big|
\frac{\beta }{N}H_N(m)+F_{N,\beta,k}(m,\delta,\epsilon) - F_{N,\beta}
\Big|<\eta
\Big)=1\,.
\end{equation}
In fact, for any $q\in[0,1]$ the limit above is either $1$ or $0$.\footnote{We do not prove this here. We only mention that, using similar arguments to those in the next two subsections, it can be shown that \eqref{eq:eta} holds if  \eqref{eq:TAPrepresentation} does and \eqref{eq:eta} holds with limit $0$ instead of $1$ if  \eqref{eq:TAPrepresentation} holds with strict inequality $>$ instead of equality.
} 
Those values of $q$ are characterized by the following.
\begin{definition}\label{def:multisamplable}
	We say that $q\in[0,1)$ is multi-samplable if for any $\epsilon>0$ and $k\geq1$,
	\begin{equation*}
		\lim_{N\to\infty}\frac1N \log \E G_{N,\beta}^{\otimes k}\Big(\forall i\neq j:\,
		\,\big|R(\bx^i,\bx^j)-q\big|<\epsilon 
		\Big)= 0\,.
	\end{equation*}
\end{definition}
\begin{lemma}\label{eq:muli-characterization}
	$q\in[0,1)$ satisfies \eqref{eq:eta}  for any $\epsilon,\delta,\eta>0$ and $k\geq1$ if and only if it is multi-samplable.
\end{lemma}
The lemma can be proved by a modification of the argument in the proof of \cite[Theorem 3]{FElandscape}. Roughly, the idea is as follows. Assume that \eqref{eq:eta} holds and let $m$ be a point as in the same equation. Then the three quantities in \eqref{eq:FkF1inequality} are within distance $\eta$ from each other.
By \eqref{eq:FbandF1}, $G_{N,\beta}(\B(m,\delta))\geq e^{-N\eta}$. The
left-hand side of \eqref{eq:orthpts} is at least $-k\eta$ since it is equal to \eqref{eq:Gk}. Hence, $q$ is multi-samplable since $k$ and $\eta$ are arbitrary and
\begin{equation}\label{eq:Gk2}
	\frac1N\log G_{N,\beta}^{\otimes k}\Big(\forall i\neq j:\,
	\,\big|R(\bx^i,\bx^j)-q\big|<\epsilon \Big)\geq -2k\eta \,.
\end{equation}

On the other hand, if $q$ is multi-samplable, then with probability not exponentially small, for independent samples $\bx^1,\ldots,\bx^{k'}$ from $G_{N,\beta}$, the pairwise overlaps satisfy $|R(\bx^i,\bx^j)-q|<\epsilon'$. Defining $m$ as as the projection to $\sqrt q\SN$ of the average $\sum_{i\leq k'}\bx^i/k'$, for $k'\gg k$ and $\epsilon\gg \epsilon'$, by simple linear algebra one can show that $(\bx^1,\ldots,\bx^{k})\in \B_{k}(m,\delta,\epsilon)$. Since $\bx^1,\ldots,\bx^{k}$ are independent samples, one also has that with probability not exponentially small such $m$ satisfies the inequality as in \eqref{eq:eta}.


\subsection{\label{sec:TAPcorrection}Generalized TAP correction.} 

We define $F_{\TAP}(q)=F_{\TAP,\beta}(q)$ as the limit\footnote{The infimum is obtained as $\epsilon,\delta\to0$ and $k\to\infty$.}
\begin{equation}\label{eq:FTAP}
F_{\TAP}(q):=\inf_{\delta,\epsilon,k} \lim_{N\to\infty}\E F_{N,\beta,k}(m,\delta,\epsilon)=\frac12\log(1-q)+F(q,\beta),
\end{equation}
where $m=m_N$ is an arbitrary sequence such that $\|m\|=\sqrt{Nq}$ and $F(q,\beta)$ is the (non-random) limiting free energy of the spherical model with mixture
\begin{equation}\label{eq:xiq}
	\xi_q(t) = \xi(q+(1-q)t)-\xi(q)-\xi'(q)(1-q)t\,.
\end{equation} 
If $q$ is either the maximal multi-samplable overlap or the rightmost point in the support of the Parisi measure, it was proved in Corollaries 6 and 12 of \cite{FElandscape} that $F(q,\beta)$ is equal to the classical Onsager correction
\begin{equation}\label{eq:Onsager}
	F(q,\beta) = \frac12\beta^2\xi_q(1)\,.
\end{equation}
This value is, of course, the free energy of the model with mixture $\xi_q(t)$ in the replica symmetric phase.
The second equality in  \eqref{eq:FTAP} was proved in \cite[Proposition 4]{FElandscape}. We now outline some of the main ideas in the proof.

For $\delta=0$, $\B(m,0)$ is a sphere centered at $m$ of dimension one less than $\SN$. By projecting $\B(m,\delta)$ to $\B(m,0)$ and mapping the latter to $\mathbb{S}^{N-2}$, we can relate $F_{N,\beta,k}(m,\delta,\epsilon)$ to a free energy on the sphere. This process incurs errors that vanish as $\delta\to0$ and we get
\begin{equation*}
	\inf_{\delta,\epsilon,k} \lim_{N\to\infty}\E F_{N,\beta,k}(m,\delta,\epsilon)=\frac12\log(1-q)+\inf_{\epsilon,k} \lim_{N\to\infty}\E \tilde F^q_{N,\beta,k}(0,1,\epsilon)\,.
\end{equation*}
Here, the logarithmic term comes from change of volumes and is equal to $$\lim\limits_{\delta\to0}\lim\limits_{N\to\infty}\frac1N\log\text{Vol}(\B(m,\delta))$$ and $\tilde F^q_{N,\beta,k}(0,1,\epsilon)$ is defined as in \eqref{eq:TAPk} with $H_N(\bx)$ replaced by the Hamiltonian $\tilde H^q_N(\bx)$ with mixture $\tilde \xi_q(t)=\xi(q+(1-q)t)-\xi(q)$. Since $m=0$ in the latter, the arbitrary choice $\delta=1$ has no effect.

For any $\bx\in\SN$, we have that
\begin{equation}\label{eq:1spin}
	\tilde H^q_N(\bx) = H^q_N(\bx) + \nabla \tilde H^q_N(0)\cdot \bx,
\end{equation}
where $H_N^q(\bx)$ is the Hamiltonian with mixture $\xi_q(t)$, independent of $\nabla \tilde H^q_N(0)$.
If $|R(\bx^i,\bx^j)|<\epsilon$ for $i<j\leq k$, 
\begin{equation}
	\label{eq:extbd}
	\Big|\frac1k\sum_{i\leq k} \nabla \tilde H^q_N(0)\cdot \bx^i\Big| 
	\leq 
	\Big\|\frac1k \sum_{i\leq k} \bx^i \Big\| \cdot\big\|\nabla \tilde H^q_N(0)\big\|\leq \sqrt{N(1/k+\epsilon)}\big\|\nabla \tilde H^q_N(0)\big\|=\sqrt{1/k+\epsilon}\cdot O(N)\,,
\end{equation}
where the equality holds with high probability. Hence, defining $F^q_{N,\beta,k}(0,1,\epsilon)$ similarly to $\tilde F^q_{N,\beta,k}(0,1,\epsilon)$  with $\tilde H^q_N(\bx)$ replaced by $H^q_N(\bx)$, 
\begin{equation}
	\label{eq:TAPext}
	\inf_{\epsilon,k} \lim_{N\to\infty}\E \tilde F^q_{N,\beta,k}(0,1,\epsilon)=\inf_{\epsilon,k} \lim_{N\to\infty}\E F^q_{N,\beta,k}(0,1,\epsilon)\,.
\end{equation}

The move from $\tilde H^q_N(\bx)$ to $H^q_N(\bx)$ in \eqref{eq:TAPext} amounts to removing the $1$-spin interaction interaction from $\tilde\xi_q(t)$ to obtain $\xi_q(t)$. For spherical models with no $1$-spin interaction one can show that the free energy  with orthogonal replicas coincides in the $N\to\infty$ limit with the usual free energy without replicas. Namely, $\lim\limits_{N\to\infty}\E F^q_{N,\beta,k}(0,1,\epsilon) = F(q,\beta)$.

\subsection{\label{sec:TAPapproach}Generalized TAP approach.} 
One of the most influential papers in the theory of spin glasses is the work  of Thouless, Anderson and Palmer (TAP) \cite{TAP}, 
where they derived the 
celebrated mean field equations for the local magnetizations of the Sherrington-Kirkpatrick model. To analyze the thermodynamics, they proposed, one should investigate the solutions of those equations, to which they associated a certain free energy $\bar F_{\rm TAP}(m)$, computed from a diagrammatic expansion of the partition function. 
Their approach was further developed in numerous works in physics, see e.g. \cite{Bray1980,Cavagna2003,CrisantiSommersTAPpspin,DeDominicis1983,Gross1984,KurchanParisiVirasoro,Plefka},  with the general idea that 
\begin{equation}
	F_{N,\beta}\approx \max_{m}\bigg\{ 
	\frac{\beta}{N} H_N(m) + \bar F_{\rm TAP}(m)
	\bigg\}, \label{eq:TAPoriginal}
\end{equation}
where the maximum is over solutions of the aforementioned equations which satisfy certain convergence conditions \cite{Plefka}. In physics, for the spherical models usually the Onsager correction $\bar F_{\rm TAP}(m)=\frac12\log(1-q)+\frac12\beta^2\xi_q(1)$ of \eqref{eq:Onsager} is used. The TAP equations are the critical point equations for the right-hand side of \eqref{eq:TAPoriginal}. The aforementioned convergence conditions define a region in $\BN$. Neglecting its boundary, equivalently to \eqref{eq:TAPoriginal} one can maximize over this region since the maximum is anyway obtained at a critical point.

Indeed, by combining \eqref{eq:FkF1inequality}, \eqref{eq:FTAP} and Theorem \ref{thm:concentration},\footnote{One also needs to prove \eqref{eq:FTAP} uniformly in $q<1-t$ for arbitrary $t>0$, but this can be done. See \cite[Proposition 4]{FElandscape}.} we obtain that w.h.p.
\begin{equation}\label{eq:myTAP}
	F_{N,\beta} =\max_{m}\bigg\{ 
	\frac{\beta}{N} H_N(m) + F_{\rm TAP}(\|m\|^2/N)
	\bigg\}+o_N(1)\,.
\end{equation}
Assuming that $\gamma_1=0$, from concentration, the origin $m=0$ is an approximate maximizer. From Lemma \ref{eq:muli-characterization}, w.h.p.
\begin{equation}\label{eq:myTAPq}
	F_{N,\beta} =\max_{m\in \sqrt{q}\SN}\bigg\{ 
	\frac{\beta}{N} H_N(m) + F_{\rm TAP}(\|m\|^2/N)
	\bigg\}+o_N(1)
\end{equation}
if and only if $q$ is multi-samplable, and otherwise we have a strict inequality $>$. We note that in \cite[Section 1.2.2]{FElandscape}, the critical point generalized TAP equations for the right-hand side of \eqref{eq:myTAP} are analyzed.

By definition, $F_{\rm TAP}(\|m\|^2/N)$ is constant on $\sqrt{q}\SN$. Define the corresponding ground-state energy,
\begin{equation*}
	\Es(q):=\lim_{N\to\infty}\frac{1}{N}\max_{\bx\in \sqrt{q}\SN}H_{N}(\bx)\,.
\end{equation*}
By taking limits in \eqref{eq:myTAPq}, we obtain the \emph{generalized TAP representation} for the free energy \cite[Theorem 5]{FElandscape}: we have that
\begin{equation}\label{eq:TAPrepresentation}
	F(\beta) = \beta \Es(q) + \frac{1}{2}\log(1-q) +F(q,\beta)
\end{equation}
if and only if $q$ is multi-samplable, while otherwise \eqref{eq:TAPrepresentation} holds as 
a strict inequality $>$. The free energy of the pure $p$-spin models was computed using this representation in \cite{SubagTAPpspin}. An analogous computation for the pure $p$-spin multi-species models was done in  \cite{Subag2023a,Subag2025}. 

Note that on $\sqrt q\SN$ with multi-samplable $q$, for large $N$,
\begin{equation} \label{eq:maxHTAP}
	\frac{\beta}{N} H_N(m) + F_{\rm TAP}(\|m\|^2/N) \approx F_{N,\beta} \iff \frac{1}{N} H_N(m) \approx \Es(q)\,.
\end{equation}
We reiterate that, by concentration of $F_{N,\beta}$ and as in Theorem \ref{thm:concentration}, points $m$ as in \eqref{eq:maxHTAP} are the ones satisfying \eqref{eq:FkF1apx}. This answers the question raised in the beginning of the section and 
provides rather simple characterizations for bands with the two geometric properties \eqref{eq:Gband} and \eqref{eq:orthpts}. In words, bands that approximately have the same free energy as the whole space and for which the probability that many independent samples from the Gibbs measure are approximately orthogonal (after centering by $m$) is not exponentially small.

Chen, Panchenko and the author \cite{TAPChenPanchenkoSubag,TAPIIChenPanchenkoSubag} proved analogues of the results in this section for mixed $p$-spin models with Ising spins. There are crucial differences, however, which make the computation of the TAP correction \eqref{eq:FTAP} (or even the proof of the existence of the limit)  significantly more challenging in the Ising case. First, because of the more complicated geometry of the cube, $\B(m,\delta)$ is generally not a simple space. Second, while similarly to \eqref{eq:TAPext} we can remove the external field ($1$-spin) term from the Hamiltonian on the band without affecting the limit of the TAP correction, in contrast to  the spherical case after its removal the free energy with multiple replicas does not have the same limit as with one replica. Overcoming these issues required new deep ideas in \cite{TAPChenPanchenkoSubag,TAPIIChenPanchenkoSubag}.

Importantly, any $q$ in the support of the Parisi measure is multi-samplable, see Corollary 8 and Theorem 10 in \cite{FElandscape}. Hence, all the results above hold for such overlap values. As we already mentioned, for the largest point in the support $F(q,\beta)$ is the classical Onsager correction \eqref{eq:Onsager}. For smaller points in the support, the free energy $F(q,\beta)$ is not replica symmetric and the corresponding Parisi measure can be expressed using that of the original model \cite[Proposition 11]{FElandscape}.

We mention earlier results, all of which concern the TAP approach with the RS classical Onsager correction for the maximal point in the support. 
For the SK model at high temperature, Talagrand \cite{TalagrandBookI} and Chatterjee \cite{ChatterjeeTAP} established the TAP equations and Bolthausen \cite{BolthausenTAP} proved a scheme to construct the solutions and used it \cite{BolthausenMorita} to compute the free energy. A generalization of the scheme is studied by Chen and Tang \cite{MR4294289,MR4718595}. An estimation problem in the setting of the SK model with a signal was studied by Fan, Mei and Montanari \cite{MR4203332}. For mixed models with Ising spins, Auffinger and Jagannath proved the TAP equations at low temperature in a formulation different from ours \cite{AuffingerJagannthSpinDist,AuffingerJagannathTAP}.
A TAP representation for the free energy for general mixed models with Ising spins was established by Chen and Panchenko \cite{ChenPanchenkoTAP}. For spherical models it was established for large $\beta$ in 
\cite{geometryMixed,geometryGibbs} which we discussed in Section \ref{sec:crtpts} and the $2$-spin model at any $\beta$ by Belius and Kistler \cite{BeliusKistlerTAP}.

\subsection{\label{sec:tree}Energies on the tree of pure states.}

In \cite{TalagrandPstates} Talagrand formulated and proved an asymptotic pure states decomposition for the Gibbs measure, which holds e.g. assuming  that $\xi(t)$ is generic \eqref{eq:generic} and the Parisi measure charges the rightmost point in its support. Jagannath \cite{JagannathApxUlt} gave another proof which removed the assumption on the Parisi measure. The decomposition states that there exists a random sequence of disjoint subsets $A_i$ such that w.h.p. $G_{N,\beta}(\cup_i A_i)=1-o(1)$ 
and the overlap distribution under the restriction of $G_{N,\beta}$ to any $A_i$  concentrates at one point.
By ultrametricity \cite{ultramet}, the pure states can be clustered according to their overlaps. This structure can be encoded through a tree with the following properties. Its vertex set $V$ is a subset of the ball $\BN$ and each pure state $A_i$ corresponds to a leaf $m_i\in V$. Denoting the least common ancestor of $m_i$ and $m_j$ by $m_i\wedge m_j$, for any $\delta>0$,
\begin{equation*}
	\lim_{N\to\infty}\E\left\{ \max_{i,j} G_{N,\beta}^{\otimes 2}\Big( \big|R(\bx^1,\bx^2)-R(m_i\wedge m_j,m_i\wedge m_j) \big|>\delta\,\Big|\,\bx^1\in A_i,\bx^2\in A_j \Big) \right\} =0.
	\label{eq:1402}
\end{equation*}	
For precise statements, we refer the reader to \cite{TalagrandPstates,JagannathApxUlt} and \cite[Corollary 14]{FElandscape}, or the work of Dembo and the author \cite[Section 3]{DemboSubagGibbs} where for $k$-RSB models the decomposition is proved with the additional property that all $\bv\in V$ are critical points.\footnote{Although the sense in which $\bv$ are critical points is slightly different in \cite{DemboSubagGibbs}, it can be easily concluded from their results that each $\bv$ is a critical point of $H_N(\bx)$ on the sphere of radius $\|\bv\|$.}
In fact, each $\bv\in V$ is a center of a band satisfying
\eqref{eq:Gband} and \eqref{eq:orthpts} for $\epsilon,\delta\ll 1$ and $k\gg1$. From the same argument as the discussion around \eqref{eq:maxHTAP}, it was proved in \cite[Corollary 15]{FElandscape} that, in probability,
\begin{equation}\label{eq:maxontree}
	\lim_{N\to\infty}\max_{\bv\in V}\Big|
	\frac1NH_N(\bv) - \max_{\|\bx\|=\|\bv\|}\frac{1}{N}H_N(\bx)\Big|=
	\lim_{N\to\infty}\max_{\bv\in V}\Big|
	\frac1NH_N(\bv) - \Es(\|\bv\|^2/N)\Big|=0\,.
\end{equation}
Another consequence of the latter corollary (cf. Theorems \ref{thm:Geometry}, \ref{thm:Geometrymixed}) is that for any multi-samplable $q$ and $\delta>0$,
\begin{equation*}
\lim_{N\to\infty} \E\Big\{ G_{N,\beta}\Big({\scalebox{1.5}{\raisebox{-0.3ex}{$\cup$}}}\Big\{\B(m,\delta):\,\|m\|=\sqrt{Nq},\,\Big|\frac1NH_N(m) - \Es(q) \Big|\leq \delta \Big\} \Big) \Big\}= 1 \,.
\end{equation*}

\section{Optimization algorithms.\label{sec:algorithms}} Suppose that for a given mixture $\xi(t)$, the support of the Parisi measure is an interval $[0,q]$ for any sufficiently large $\beta$. Such mixtures are characterized by the condition that $\xi''(t)^{-1/2}$ is concave on $(0,1]$, see \cite[Prop. 1]{GSFollowing}. In this case, the tree described in Subsection \ref{sec:tree} can be constructed so that w.h.p. any root-to-leaf path $\bv_0,\bv_1,\ldots,\bv_{k(N)}$ satisfies the following, where $k(N)$ is the depth of the tree which diverges with $N$. First, it has the geometric properties that $R(\bv_i,\bv_i)=iq/k(N)$ for any $i$ and $R(\bv_j-\bv_i,\bv_n-\bv_l)=o(1)$ for any $i\leq j\leq l\leq n$. Second, due to \eqref{eq:maxontree}, it approximately maximizes the energy (consistently) in the sense that 
\begin{equation}
	\big|H_N(\bv_i)/N - \Es(iq/k(N))\big|=o(1)\,.
\end{equation} 
Since as $\beta\to\infty$ the rightmost point in the support goes to $1$, such paths exist also for $q=1$.

This geometric structure suggests a natural optimization algorithm: to find a point approximately maximizing $H_N(\bx)$ on $\SN$, construct a path as above and output its end point.
The geometric properties of the path can be  imposed as follows. Fixing some $k\geq1$ and setting $\delta=1/k$,  define a sequence $\bx_0,\ldots\bx_k$ such that the increments 
$\bu_i:=\bx_{i+1}-\bx_{i}$ satisfy $\|\bu_i\|=\sqrt{N\delta}$ and $\bu_i\perp \bx_i$, so that $\|\bx_i\|^2=i\delta N$.\footnote{Here we only require orthogonality of the increment and current position. We could in fact impose orthogonality between all increments as for the path on the ultrametric tree, but this will not be important for the algorithm.}  In particular, the output $\bxs:=\bx_k$ belongs to $\SN$. The directions $\hat\bu_i = \bu_i/\|\bu_i\|$ need to be chosen from  $\bx_i^\perp$,  the orthogonal space to $\bx_i$, in a way that keeps the energy approximately maximal on the corresponding sphere.  
It may seem that a natural choice is to take the direction of the gradient  $\nabla H_N(\bx_i)$ projected to $\bx_i^\perp$. However, if our goal of having $H_N(\bx_i)$ approximately maximal over the sphere of radius $\|\bx_i\|$ is achieved, this projected gradient should be close to $0$ (after scaling by $\sqrt N$). We may instead choose a direction that maximizes the energy increment at second order. Namely, choose $\hat\bu_i$ such that, for small $\epsilon>0$,
\begin{equation*}
	\hat\bu_i^T \nabla^2 H_N(\bx_i) \hat\bu_i \geq \max_{\bu\perp \bx_i, \|\bu\|=1} \bu^T \nabla^2 H_N(\bx_i) \bu - \epsilon\,.
\end{equation*}  

It was proved  in \cite{GSFollowing} that, under the assumption above on the support of the Parisi measure, this algorithm indeed produces $\bxs\in\SN$ such that $H_N(\bxs)/N \geq \Es-\eta$, where $\eta\to0$ as $\epsilon,\delta\to0$. It is not difficult to check that it has polynomial in $N$ running time. For general mixtures $\xi(t)$, the same algorithm achieves
\begin{equation}
	H_N(\bxs) \geq \int_0^1 \sqrt{\xi''(t)}dt-\eta\,.
\end{equation}

Recent years have seen significant advances in optimization of spin glasses, with several noteworthy developments. For models with Ising spins, Montanari \cite{MontanariSKopt} introduced an approximate-message-passing algorithm that optimizes the energy of the SK model and  El Alaoui, Montanari and Sellke \cite{ElAlaoui2021} developed and studied a variant of the latter for general mixtures. In particular, they derived a Parisi-like formula for the terminal energy of these algorithms. Like the Hessian descent algorithm described above for the spherical models, the algorithms of \cite{ElAlaoui2021,MontanariSKopt} also construct a path with orthogonal increments from the origin to the configuration space\,---\,the discrete cube, in the Ising case\,---\,a fact that was crucial to the analysis of those algorithms.
Huang and Sellke \cite{HuangSellkeHardness} proved that within a class of algorithms that are Lipschitz in an appropriate sense, the algorithms of \cite{GSFollowing,MontanariSKopt,ElAlaoui2021} achieve the best energy. See also \cite{HuangSellkeHardnessMulti,MR4709433} for their related works on multi-species spherical models. It was conjectured in \cite{GSFollowing} that applying a Hessian descent algorithm similar to the above to the modified energy $\frac1NH_N(\bx)+F_\TAP(\bx)$ (where $F_\TAP(\bx)$ is the generalized TAP correction), should result in an algorithm that approximately achieves the maximal energy for full-RSB models also in the Ising case. For the SK model, Jekel, Sandhu and Shi recently resolved this conjecture in \cite{Jekel2025}.

\subsection{Polynomial systems of equations.}
In 1998 Steve Smale published a list of ‘Mathematical problems for the next century’ \cite{SmaleProblems}. His 17th problem asked for a poly-time algorithm for solving random polynomial systems. Research on the problem has focused on the complex case, which was resolved after several breakthroughs in 2017, see the works of Beltr\'{a}n and Pardo \cite{Beltran2008}, B\"{u}rgisser and Cucker \cite{Buergisser2011} and Lairez \cite{Lairez2017} and references therein. No significant progress was made, however, on the more difficult case over the real numbers, to which previous methods cannot be applied.
In the real case, the polynomial system in Smale's problem consists of $N-1$ independent random homogeneous polynomials in $N$ variables, where each polynomial follows the law of a pure $p$-spin Hamiltonian normalized by $\sqrt N$ (with $p$ possibly varying across equations). The input to the algorithm is the set of coefficients of all the polynomials, and the output needs to be a point from which Newton's method on $\SN$ converges quadratically to an exact solution. The model of computation is over the real numbers\,---\,formally, a Blum-Shub-Smale machine. 

Denoting the system by $\bF(\bx)=(F_1(\bx),\ldots,F_{N-1}(\bx))$ and defining an energy function $H(\bx):=\|\bF(\bx)\|^2$, note that $H(\bx)=0$ if and only if $\bF(\bx)=\mathbf{0}$. I.e., solving $\bF(\bx)$ is equivalent to finding a global minimum point of $H(\bx)$. Montanari and the author \cite{MontanariSubag_Smale17} introduced a variant of the Hessian descent optimization algorithm discussed above and resolved a simplified version of Smale's 17th problem with $N-O(\sqrt{N\log N})$ equations instead of $N-1$. 

In another work \cite{MontanariSubag2023}, the same authors study the problem of minimizing $H(\bx)/N$ over $\SN$ in the $N\to\infty$ limit for $H(\bx) = \sum_{i=1}^{\lfloor\alpha N\rfloor} |F_i(\bx)-c|^2$, with $F_i(\bx)$ being i.i.d. copies of a mixed $p$-spin Hamiltonian with mixture $\xi(t)$ and parameters $\xi,c,\alpha$ held fixed. This corresponds to approximately solving the system $F_i(\bx)=c$.\footnote{To be precise, in \cite{MontanariSubag2023} the authors consider the case $c=0$, but allow the mixed models include a $0$-spin (constant) term, which shift each $F_i(\bx)$ by an independent  Gaussian variable that is uniform over the sphere $\SN$. Due to rotational invariance, asymptotically this is equivalent to our setting without the $0$-spin term and with fixed $c\geq 0$.} Note that for sufficiently large $c>0$, the system $F_i(\bx)=c$ has no solution with high probability. In this setting, it is proved in \cite{MontanariSubag2023} that a Hessian descent minimization algorithm produces a path $\bx_i$ from the origin to the sphere such that $H_N(\bx_i)/N\leq \varphi(\|\bx_i\|^2/N)+o(1)$, where $\varphi(t)$ solves an explicit differential equation. If $\varphi(1)=0$, then the terminal point approximately solves $F_i(\bx)=c$.


\bibliographystyle{siamplain}
\bibliography{master}

\begin{thebibliography}{100}

\bibitem{RFG}
{\sc R.~J. Adler and J.~E. Taylor}, {\em Random fields and geometry}, Springer
  Monographs in Mathematics, Springer, New York, 2007.

\bibitem{AizenmanContucci}
{\sc M.~Aizenman and P.~Contucci}, {\em On the stability of the quenched state
  in mean-field spin-glass models}, J. Statist. Phys., 92 (1998), pp.~765--783.

\bibitem{AS2scheme}
{\sc M.~Aizenman, R.~Sims, and S.~L. Starr}, {\em Extended variational
  principle for the {S}herrington-{K}irkpatrick spin-glass model}, Phys. Rev.
  B, 68 (2003), p.~214403.

\bibitem{Anderson1989}
{\sc P.~W. Anderson}, {\em {Spin Glass VI: Spin Glass As Cornucopia}}, Physics
  Today, 42 (1989), pp.~9--11.

\bibitem{Montanari_emprisk}
{\sc K.~Asgari, A.~Montanari, and B.~Saeed}, {\em Local minima of the empirical
  risk in high dimension: General theorems and convex examples},
  arXiv:2502.01953,  (2025).

\bibitem{ABA2}
{\sc A.~Auffinger and G.~Ben~Arous}, {\em Complexity of random smooth functions
  on the high-dimensional sphere}, Ann. Probab., 41 (2013), pp.~4214--4247.

\bibitem{A-BA-C}
{\sc A.~Auffinger, G.~Ben~Arous, and J.~{\v{C}}ern{\'y}}, {\em Random matrices
  and complexity of spin glasses}, Comm. Pure Appl. Math., 66 (2013),
  pp.~165--201.

\bibitem{MR4404771}
{\sc A.~Auffinger, G.~Ben~Arous, and Z.~Li}, {\em Sharp complexity asymptotics
  and topological trivialization for the {$(p,k)$} spiked tensor model}, J.
  Math. Phys., 63 (2022), pp.~Paper No. 043303, 21.

\bibitem{MR3249903}
{\sc A.~Auffinger and W.-K. Chen}, {\em Free energy and complexity of spherical
  bipartite models}, J. Stat. Phys., 157 (2014), pp.~40--59.

\bibitem{AuffingerChenConcentration}
{\sc A.~Auffinger and W.-K. Chen}, {\em On concentration properties of
  disordered hamiltonians}, Proc. Am. Math. Soc, 146 (2018), pp.~1807--1815.

\bibitem{AuffingerJagannthSpinDist}
{\sc A.~Auffinger and A.~Jagannath}, {\em On spin distributions for generic
  p-spin models}, J. Stat. Phys., 174 (2019), pp.~316--332.

\bibitem{AuffingerJagannathTAP}
{\sc A.~Auffinger and A.~Jagannath}, {\em {T}houless-{A}nderson-{P}almer
  equations for generic $p$-spin glasses}, Ann. Probab., 47 (2019),
  pp.~2230--2256.

\bibitem{Azais2005}
{\sc J.-M. Aza\"{\i}s and M.~Wschebor}, {\em On the roots of a random system of
  equations. {T}he theorem of {S}hub and {S}male and some extensions}, Found.
  Comput. Math., 5 (2005), pp.~125--144.

\bibitem{BaffioniRosati}
{\sc F.~Baffioni and F.~Rosati}, {\em Some exact results on the ultrametric
  overlap distribution in mean field spin glass models (i)}, Eur. Phys. J. B,
  17 (2000), pp.~439--447.

\bibitem{BeliusTAP1}
{\sc D.~Belius, F.~Concetti, and G.~Genovese}, {\em On the determinant in
  {Bray-Moore's TAP} complexity formula}, arXiv:2401.08529,  (2024).

\bibitem{BeliusKistlerTAP}
{\sc D.~Belius and N.~Kistler}, {\em The {TAP}-{P}lefka variational principle
  for the spherical {SK} model}, Comm. Math. Phys., 367 (2019), pp.~991--1017.

\bibitem{BeliusSchmidt}
{\sc D.~Belius and M.~A. Schmidt}, {\em Complexity of local maxima of given
  radial derivative for mixed $p$-spin {H}amiltonians}, arXiv:2207.14361,
  (2022).

\bibitem{MR4354698}
{\sc D.~Belius, J.~\v{C}ern\'{y}, S.~Nakajima, and M.~A. Schmidt}, {\em
  Triviality of the geometry of mixed {$p$}-spin spherical {H}amiltonians with
  external field}, J. Stat. Phys., 186 (2022), pp.~Paper No. 12, 34.

\bibitem{Beltran2008}
{\sc C.~Beltr\'{a}n and L.~M. Pardo}, {\em On {S}male's 17th problem: a
  probabilistic positive solution}, Found. Comput. Math., 8 (2008), pp.~1--43,
  \url{https://doi.org/10.1007/s10208-005-0211-0},
  \url{https://doi.org/10.1007/s10208-005-0211-0}.

\bibitem{MR4552227}
{\sc G.~Ben~Arous, P.~Bourgade, and B.~McKenna}, {\em Exponential growth of
  random determinants beyond invariance}, Probab. Math. Phys., 3 (2022),
  pp.~731--789.

\bibitem{MR4673883}
{\sc G.~Ben~Arous, P.~Bourgade, and B.~McKenna}, {\em Landscape complexity
  beyond invariance and the elastic manifold}, Comm. Pure Appl. Math., 77
  (2024), pp.~1302--1352.

\bibitem{MR4305622}
{\sc G.~Ben~Arous, Y.~V. Fyodorov, and B.~A. Khoruzhenko}, {\em Counting
  equilibria of large complex systems by instability index}, Proc. Natl. Acad.
  Sci. USA, 118 (2021).

\bibitem{BenArousJagannathShattering}
{\sc G.~Ben~Arous and A.~Jagannath}, {\em Shattering versus metastability in
  spin glasses}, Comm. Pure Appl. Math., 77 (2024), pp.~139--176.

\bibitem{MR4011861}
{\sc G.~Ben~Arous, S.~Mei, A.~Montanari, and M.~Nica}, {\em The landscape of
  the spiked tensor model}, Comm. Pure Appl. Math., 72 (2019), pp.~2282--2330.

\bibitem{geometryMixed}
{\sc G.~Ben~Arous, E.~Subag, and O.~Zeitouni}, {\em Geometry and temperature
  chaos in mixed spherical spin glasses at low temperature: the perturbative
  regime}, Comm. Pure Appl. Math., 73 (2020), pp.~1732--1828.

\bibitem{Blandin}
{\sc A.~Blandin}, {\em Theories versus experiments in spin glass systems}, J.
  Phys. Colloq., 39 (1978).

\bibitem{BolthausenTAP}
{\sc E.~Bolthausen}, {\em An iterative construction of solutions of the {TAP}
  equations for the {S}herrington-{K}irkpatrick model}, Comm. Math. Phys., 325
  (2014), pp.~333--366.

\bibitem{BolthausenMorita}
{\sc E.~Bolthausen}, {\em A {M}orita type proof of the replica-symmetric
  formula for {SK}}, in Statistical Mechanics of Classical and Disordered
  Systems, vol.~293 of Springer Proceedings in Mathematics \& Statistics,
  Springer, 2019.

\bibitem{BrayMoorePhysRevLett.41.1068}
{\sc A.~J. Bray and M.~A. Moore}, {\em Replica-symmetry breaking in spin-glass
  theories}, Phys. Rev. Lett., 41 (1978), pp.~1068--1072.

\bibitem{Bray1980}
{\sc A.~J. Bray and M.~A. Moore}, {\em Metastable states in spin glasses}, J.
  Phys. C, 13 (1980), p.~L469.

\bibitem{Buergisser2011}
{\sc P.~B\"{u}rgisser and F.~Cucker}, {\em On a problem posed by {S}teve
  {S}male}, Ann. of Math. (2), 174 (2011), pp.~1785--1836.

\bibitem{Cavagna2003}
{\sc A.~Cavagna, I.~Giardina, G.~Parisi, and M.~M{\'e}zard}, {\em On the formal
  equivalence of the {TAP} and thermodynamic methods in the {SK} model}, J.
  Phys. A, 36 (2003), pp.~1175--1194.

\bibitem{CharbonneauReplicaReview}
{\sc P.~Charbonneau}, {\em From the replica trick to the replica symmetry
  breaking technique}, arXiv:2211.01802,  (2022).

\bibitem{ChatterjeeTAP}
{\sc S.~Chatterjee}, {\em Spin glasses and {S}tein's method}, Probab. Theory
  Related Fields, 148 (2010), pp.~567--600.

\bibitem{Chen2013}
{\sc W.-K. Chen}, {\em The {A}izenman-{S}ims-{S}tarr scheme and {P}arisi
  formula for mixed {$p$}-spin spherical models}, Electron. J. Probab., 18
  (2013), pp.~no. 94, 14.

\bibitem{ChenPanchenkoTAP}
{\sc W.-K. Chen and D.~Panchenko}, {\em On the {TAP} free energy in the mixed
  {$p$}-spin models}, Comm. Math. Phys., 362 (2018), pp.~219--252.

\bibitem{TAPChenPanchenkoSubag}
{\sc W.-K. Chen, D.~Panchenko, and E.~Subag}, {\em The generalized {TAP} free
  energy}, Comm. Pure Appl. Math., 76, pp.~1329--1415.

\bibitem{TAPIIChenPanchenkoSubag}
{\sc W.-K. Chen, D.~Panchenko, and E.~Subag}, {\em The generalized {TAP} free
  energy {II}}, Comm. Math. Phys., 381 (2021), pp.~257--291.

\bibitem{MR4294289}
{\sc W.-K. Chen and S.~Tang}, {\em On convergence of the cavity and
  {B}olthausen's {TAP} iterations to the local magnetization}, Comm. Math.
  Phys., 386 (2021), pp.~1209--1242.

\bibitem{MR4718595}
{\sc W.-K. Chen and S.~Tang}, {\em On the {TAP} equations via the cavity
  approach in the generic mixed {$p$}-spin models}, Comm. Math. Phys., 405
  (2024), pp.~Paper No. 87, 43.

\bibitem{Crisanti1992}
{\sc A.~Crisanti and H.-J. Sommers}, {\em The spherical p-spin interaction spin
  glass model: the statics}, Z. Physik B - Condensed Matter, 87 (1992),
  pp.~341--354.

\bibitem{CrisantiSommersTAPpspin}
{\sc A.~Crisanti and H.-J. Sommers}, {\em Thouless-anderson-palmer approach to
  the spherical p-spin spin glass model}, J. Phys. I France, 5 (1995),
  pp.~805--813.

\bibitem{DartoisMcKenna}
{\sc S.~Dartois and B.~McKenna}, {\em Injective norm of real and complex random
  tensors {I}: From spin glasses to geometric entanglement}, arXiv:2404.03627,
  (2024).

\bibitem{Almeida1978}
{\sc J.~R.~L. de~Almeida and D.~J. Thouless}, {\em Stability of the
  {Sherrington-Kirkpatrick} solution of a spin glass model}, J. Phys. A, 11
  (1978), p.~983.

\bibitem{DeDominicis1983}
{\sc C.~De~Dominicis and A.~P. Young}, {\em Weighted averages and order
  parameters for the infinite range {I}sing spin glass}, J. Phys. A, 16 (1983),
  pp.~2063--2075.

\bibitem{DemboSubagGibbs}
{\sc A.~Dembo and E.~Subag}, {\em {Disordered Gibbs measures and Gaussian
  conditioning}}, arXiv:2409.19453,  (2024).

\bibitem{Dovbysh-Sudakov}
{\sc L.~N. Dovbysh and V.~N. Sudakov}, {\em Gram-de {F}inetti matrices}, Zap.
  Nauchn. Sem. Leningrad. Otdel. Mat. Inst. Steklov. (LOMI), 119 (1982),
  pp.~77--86, 238, 244--245.
\newblock Problems of the theory of probability distribution, VII.

\bibitem{EA1975}
{\sc S.~F. Edwards and P.~W. Anderson}, {\em Theory of spin glasses}, J. Phys.
  F, 5 (1975), p.~965.

\bibitem{ElAlaoui2021}
{\sc A.~El~Alaoui, A.~Montanari, and M.~Sellke}, {\em {Optimization of
  mean-field spin glasses}}, Ann. Probab., 49 (2021), pp.~2922 -- 2960.

\bibitem{MR4203332}
{\sc Z.~Fan, S.~Mei, and A.~Montanari}, {\em T{AP} free energy, spin glasses
  and variational inference}, Ann. Probab., 49 (2021), pp.~1--45.

\bibitem{Fyodorov2004}
{\sc Y.~V. Fyodorov}, {\em Complexity of random energy landscapes, glass
  transition, and absolute value of the spectral determinant of random
  matrices}, Phys. Rev. Lett., 92 (2004), pp.~240601, 4.

\bibitem{Fyodorov2013}
{\sc Y.~V. Fyodorov}, {\em High-dimensional random fields and random matrix
  theory}, Markov Process. Relat. Fields, 21 (2015), pp.~483--518.

\bibitem{MR3598251}
{\sc Y.~V. Fyodorov}, {\em Topology trivialization transition in random
  non-gradient autonomous {ODE}s on a sphere}, J. Stat. Mech. Theory Exp.,
  (2016), pp.~124003, 21.

\bibitem{MR3521630}
{\sc Y.~V. Fyodorov and B.~A. Khoruzhenko}, {\em Nonlinear analogue of the
  {M}ay-{W}igner instability transition}, Proc. Natl. Acad. Sci. USA, 113
  (2016), pp.~6827--6832.

\bibitem{MR3162549}
{\sc Y.~V. Fyodorov and P.~Le~Doussal}, {\em Topology trivialization and large
  deviations for the minimum in the simplest random optimization}, J. Stat.
  Phys., 154 (2014), pp.~466--490.

\bibitem{MR2363390}
{\sc Y.~V. Fyodorov and I.~Williams}, {\em Replica symmetry breaking condition
  exposed by random matrix calculation of landscape complexity}, J. Stat.
  Phys., 129 (2007), pp.~1081--1116.

\bibitem{Garciacomplexity}
{\sc X.~Garcia}, {\em On the number of equilibria with a given number of
  unstable directions}, arXiv:1709.04021,  (2017).

\bibitem{GhirlandaGuerra}
{\sc S.~Ghirlanda and F.~Guerra}, {\em General properties of overlap
  probability distributions in disordered spin systems. {T}owards {P}arisi
  ultrametricity}, J. Phys. A, 31 (1998), pp.~9149--9155.

\bibitem{Gross1984}
{\sc D.~J. Gross and M.~Mezard}, {\em The simplest spin glass}, Nuclear Physics
  B, 240 (1984), pp.~431 -- 452.

\bibitem{GuerraBound}
{\sc F.~Guerra}, {\em Broken replica symmetry bounds in the mean field spin
  glass model}, Comm. Math. Phys., 233 (2003), pp.~1--12.

\bibitem{GuerraToninelli}
{\sc F.~Guerra and F.~L. Toninelli}, {\em The thermodynamic limit in mean field
  spin glass models}, Comm. Math. Phys., 230 (2002), pp.~71--79.

\bibitem{HuangSellkeHardnessMulti}
{\sc B.~Huang and M.~Sellke}, {\em Algorithmic threshold for multi-species
  spherical spin glasses}, arXiv:2303.12172,  (2023).

\bibitem{HuangSellke2023}
{\sc B.~Huang and M.~Sellke}, {\em A constructive proof of the spherical
  {P}arisi formula}, Preprint arXiv:2311.15495,  (2023).

\bibitem{Huang_toptriv}
{\sc B.~Huang and M.~Sellke}, {\em Strong topological trivialization of
  multi-species spherical spin glasses}, arXiv:2308.09677,  (2023).

\bibitem{MR4709433}
{\sc B.~Huang and M.~Sellke}, {\em Optimization algorithms for multi-species
  spherical spin glasses}, J. Stat. Phys., 191 (2024), pp.~Paper No. 29, 42.

\bibitem{HuangSellkeHardness}
{\sc B.~Huang and M.~Sellke}, {\em Tight {L}ipschitz hardness for optimizing
  mean field spin glasses}, Comm. Pure Appl. Math., 78 (2025), pp.~60--119.

\bibitem{JagannathApxUlt}
{\sc A.~Jagannath}, {\em Approximate ultrametricity for random measures and
  applications to spin glasses}, Comm. Pure Appl. Math., 70 (2017),
  pp.~611--664.

\bibitem{Jekel2025}
{\sc D.~Jekel, J.~S. Sandhu, and J.~Shi}, {\em Potential {H}essian ascent: the
  {S}herrington-{K}irkpatrick model}, in SODA, SIAM, Philadelphia, PA, 2025,
  pp.~5307--5387.

\bibitem{Kac1943}
{\sc M.~Kac}, {\em On the average number of real roots of a random algebraic
  equation}, Bull. Amer. Math. Soc., 49 (1943), pp.~314--320.

\bibitem{Kirkpatrick1983}
{\sc S.~Kirkpatrick, C.~D. Gelatt, and M.~P. Vecchi}, {\em Optimization by
  simulated annealing}, Science, 220 (1983), pp.~671--680.

\bibitem{KirkpatrickPhysRevB.17.4384}
{\sc S.~Kirkpatrick and D.~Sherrington}, {\em Infinite-ranged models of
  spin-glasses}, Phys. Rev. B, 17 (1978), pp.~4384--4403.

\bibitem{MR4645713}
{\sc P.~Kivimae}, {\em The ground state energy and concentration of complexity
  in spherical bipartite models}, Comm. Math. Phys., 403 (2023), pp.~37--81.

\bibitem{MR4824712}
{\sc P.~Kivimae}, {\em Concentration of equilibria and relative instability in
  disordered non-relaxational dynamics}, Comm. Math. Phys., 405 (2024),
  pp.~Paper No. 289, 66.

\bibitem{KurchanParisiVirasoro}
{\sc J.~Kurchan, G.~Parisi, and M.~A. Virasoro}, {\em Barriers and metastable
  states as saddle points in the replica approach}, J. Phys. I France, 3
  (1993), pp.~1819--1838.

\bibitem{MR4411381}
{\sc B.~Lacroix-A-Chez-Toine and Y.~V. Fyodorov}, {\em Counting equilibria in a
  random non-gradient dynamics with heterogeneous relaxation rates}, J. Phys.
  A, 55 (2022), pp.~Paper No. 144001, 39.

\bibitem{Lairez2017}
{\sc P.~Lairez}, {\em A deterministic algorithm to compute approximate roots of
  polynomial systems in polynomial average time}, Found. Comput. Math., 17
  (2017), pp.~1265--1292.

\bibitem{pmlr-v107-maillard20a}
{\sc A.~Maillard, G.~Ben~Arous, and G.~Biroli}, {\em Landscape complexity for
  the empirical risk of generalized linear models}, in Proceedings of The First
  Mathematical and Scientific Machine Learning Conference, vol.~107, PMLR,
  2020, pp.~287--327.

\bibitem{MR4718393}
{\sc B.~McKenna}, {\em Complexity of bipartite spherical spin glasses}, Ann.
  Inst. Henri Poincar\'e{} Probab. Stat., 60 (2024), pp.~636--657.

\bibitem{MPSTV2}
{\sc M.~M\'ezard, G.~Parisi, N.~Sourlas, G.~Toulouse, and M.~Virasoro}, {\em
  Nature of the spin-glass phase}, Phys. Rev. Lett., 52 (1984), pp.~1156--1159.

\bibitem{MPSTV1}
{\sc M.~M{\'e}zard, G.~Parisi, N.~Sourlas, G.~Toulouse, and M.~Virasoro}, {\em
  Replica symmetry breaking and the nature of the spin glass phase}, J.
  Physique, 45 (1984), pp.~843--854.

\bibitem{MPVspinglass}
{\sc M.~M{\'e}zard, G.~Parisi, and M.~A. Virasoro}, {\em Spin Glass Theory and
  Beyond}, World Scientific, 1987.

\bibitem{MontanariSKopt}
{\sc A.~Montanari}, {\em Optimization of the {S}herrington--{K}irkpatrick
  {H}amiltonian}, SIAM J. Comput., 54 (2025), pp.~1--38.

\bibitem{MontanariSubag2023}
{\sc A.~Montanari and E.~Subag}, {\em Solving systems of random equations via
  first and second-order optimization algorithms}, Preprint arXiv:2306.13326,
  (2023).

\bibitem{MontanariSubag_Smale17}
{\sc A.~Montanari and E.~Subag}, {\em {On Smale's 17th problem over the
  reals}}, arXiv:2405.01735,  (2024).

\bibitem{MR2399285}
{\sc D.~Panchenko}, {\em On differentiability of the {P}arisi formula},
  Electron. Commun. Probab., 13 (2008), pp.~241--247.

\bibitem{ultramet}
{\sc D.~Panchenko}, {\em The {P}arisi ultrametricity conjecture}, Ann. of Math.
  (2), 177 (2013), pp.~383--393.

\bibitem{PanchenkoBook}
{\sc D.~Panchenko}, {\em The {S}herrington-{K}irkpatrick model}, Springer
  Monographs in Mathematics, Springer, 2013.

\bibitem{Panchenko2014}
{\sc D.~Panchenko}, {\em The {P}arisi formula for mixed {$p$}-spin models},
  Ann. Probab., 42 (2014), pp.~946--958.

\bibitem{MR4680407}
{\sc D.~Panchenko}, {\em Ultrametricity in spin glasses}, in
  I{CM}---{I}nternational {C}ongress of {M}athematicians. {V}ol. 6. {S}ections
  12--14, EMS Press, Berlin, 2023, pp.~4376--4392.

\bibitem{ParisiFormula}
{\sc G.~Parisi}, {\em Infinite number of order parameters for spin-glasses},
  Phys. Rev. Lett., 43 (1979), pp.~1754--1756.

\bibitem{Parisi1980}
{\sc G.~Parisi}, {\em A sequence of approximated solutions to the {S-K} model
  for spin glasses}, J. Phys. A, 13 (1980), p.~L115.

\bibitem{ParisiOrderPar}
{\sc G.~Parisi}, {\em Order parameter for spin-glasses}, Phys. Rev. Lett., 50
  (1983), pp.~1946--1948.

\bibitem{ParisiNobelLecture}
{\sc G.~Parisi}, {\em Nobel lecture: Multiple equilibria}, Rev. Mod. Phys., 95
  (2023), p.~030501.

\bibitem{Plefka}
{\sc T.~Plefka}, {\em Convergence condition of the {TAP} equation for the
  infinite-ranged {I}sing spin glass model}, J. Phys. A, 15 (1982),
  pp.~1971--1978.

\bibitem{Rice1945}
{\sc S.~O. Rice}, {\em Mathematical analysis of random noise}, Bell System
  Tech. J., 24 (1945), pp.~46--156.

\bibitem{PhysRevX.9.011003}
{\sc V.~Ros, G.~Ben~Arous, G.~Biroli, and C.~Cammarota}, {\em Complex energy
  landscapes in spiked-tensor and simple glassy models: Ruggedness,
  arrangements of local minima, and phase transitions}, Phys. Rev. X, 9 (2019).

\bibitem{Ruelle}
{\sc D.~Ruelle}, {\em A mathematical reformulation of {D}errida's {REM} and
  {GREM}}, Comm. Math. Phys., 108 (1987), pp.~225--239.

\bibitem{Schoenberg}
{\sc I.~J. Schoenberg}, {\em Positive definite functions on spheres}, Duke
  Math. J., 9 (1942), pp.~96--108.

\bibitem{SK75}
{\sc D.~Sherrington and S.~Kirkpatrick}, {\em Solvable model of a spin glass},
  Phys. Rev. Lett., 35 (1975), pp.~1792--1795.

\bibitem{SmaleProblems}
{\sc S.~Smale}, {\em Mathematical problems for the next century}, Math.
  Intelligencer, 20 (1998), pp.~7--15.

\bibitem{2nd}
{\sc E.~Subag}, {\em The complexity of spherical {$p$}-spin models---{A} second
  moment approach}, Ann. Probab., 45 (2017), pp.~3385--3450.

\bibitem{geometryGibbs}
{\sc E.~Subag}, {\em The geometry of the {G}ibbs measure of pure spherical spin
  glasses}, Invent. Math., 210 (2017), pp.~135--209.

\bibitem{GSFollowing}
{\sc E.~Subag}, {\em Following the ground-states of full-{RSB} spherical spin
  glasses}, Comm. Pure Appl. Math., 74 (2021), pp.~1021--1044.

\bibitem{FreeEnergyConvergence}
{\sc E.~Subag}, {\em Convergence of the free energy for spherical spin
  glasses}, J. Stat. Phys., 189(2) (2022).

\bibitem{Subag2023}
{\sc E.~Subag}, {\em Concentration for the zero set of random polynomial
  systems}, arXiv:2303.11924,  (2023).

\bibitem{SubagTAPpspin}
{\sc E.~Subag}, {\em The free energy of spherical pure {$p$}-spin models:
  computation from the {TAP} approach}, Probab. Theory Related Fields, 186
  (2023), pp.~715--734.

\bibitem{Subag2023a}
{\sc E.~Subag}, {\em T{AP} approach for multispecies spherical spin glasses
  {II}: {T}he free energy of the pure models}, Ann. Probab., 51 (2023),
  pp.~1004--1024.

\bibitem{FElandscape}
{\sc E.~Subag}, {\em Free energy landscapes in spherical spin glasses}, Duke
  Math. J., 173 (2024), pp.~1291--1357.

\bibitem{Subag2025}
{\sc E.~Subag}, {\em {TAP} approach for multi-species spherical spin glasses
  {I}: general theory}, Electron. J. Probab., 30 (2025), pp.~Paper No. 87, 32.

\bibitem{pspinext}
{\sc E.~Subag and O.~Zeitouni}, {\em The extremal process of critical points of
  the pure {$p$}-spin spherical spin glass model}, Probab. Theory Related
  Fields, 168 (2017), pp.~773--820.

\bibitem{2ndarbitraryenergy}
{\sc E.~Subag and O.~Zeitouni}, {\em Concentration of the complexity of
  spherical pure {$p$}-spin models at arbitrary energies}, J. Math. Phys., 62
  (2021).

\bibitem{MR1648045}
{\sc M.~Talagrand}, {\em Huge random structures and mean field models for spin
  glasses}, in Proceedings of the {I}nternational {C}ongress of
  {M}athematicians, {V}ol. {I} ({B}erlin, 1998), 1998, pp.~507--536.

\bibitem{TalagrandOldBook}
{\sc M.~Talagrand}, {\em Spin glasses: a challenge for mathematicians},
  Springer-Verlag, Berlin, 2003.

\bibitem{Talag}
{\sc M.~Talagrand}, {\em Free energy of the spherical mean field model},
  Probab. Theory Relat. Fields, 134 (2006), pp.~339--382.

\bibitem{Talag2}
{\sc M.~Talagrand}, {\em The {P}arisi formula}, Ann. of Math. (2), 163 (2006),
  pp.~221--263.

\bibitem{MR2195333}
{\sc M.~Talagrand}, {\em Parisi measures}, J. Funct. Anal., 231 (2006),
  pp.~269--286.

\bibitem{TalagrandPstates}
{\sc M.~Talagrand}, {\em Construction of pure states in mean field models for
  spin glasses}, Probab. Theory Related Fields, 148 (2010), pp.~601--643.

\bibitem{TalagrandBookI}
{\sc M.~Talagrand}, {\em Mean field models for spin glasses. {V}olume {I}},
  Springer-Verlag, Berlin, 2011.

\bibitem{TalagrandBookII}
{\sc M.~Talagrand}, {\em Mean field models for spin glasses. {V}olume {II}},
  Springer, Heidelberg, 2011.

\bibitem{TAP}
{\sc D.~J. Thouless, P.~W. Anderson, and R.~G. Palmer}, {\em Solution of
  `solvable model of a spin glass'}, Philos. Mag., 35 (1977), pp.~593--601.

\bibitem{Wschebor2005}
{\sc M.~Wschebor}, {\em On the {K}ostlan-{S}hub-{S}male model for random
  polynomial systems. {V}ariance of the number of roots}, J. Complexity, 21
  (2005), pp.~773--789.

\end{thebibliography}

\end{document}